\documentclass{ws-ijbc}

\usepackage{ws-rotating} 

\usepackage{amsfonts, amsbsy, amssymb, amsmath, commath, graphicx, float, footnote, natbib, subfigure, color, url, titlesec, textcomp, gensymb, epstopdf, ifpdf, verbatim}

\usepackage[english]{babel}

\begin{document}

\title{A Theoretical Framework for Lagrangian Descriptors}
\author{C. Lopesino$^{1}$, F. Balibrea-Iniesta$^{1}$, V. J. Garc\'{i}a-Garrido$^{1}$,\\
S. Wiggins$^2$, A. M. Mancho$^1$ \vspace{0.3cm} \\
$^1$Instituto de Ciencias Matem\'aticas, CSIC-UAM-UC3M-UCM,
\\ C/Nicol\'as Cabrera 15, Campus Cantoblanco UAM, 28049
Madrid, Spain\\
$^2$School of Mathematics, University of Bristol, \\Bristol BS8 1TW, United Kingdom}
\maketitle

\begin{abstract}
This  paper provides  a theoretical background for Lagrangian Descriptors (LDs). The goal of achieving  rigourous proofs that justify the ability of LDs to detect invariant manifolds is simplified by introducing an alternative definition for LDs. The definition is stated for $n$-dimensional systems with general time dependence, however we rigorously prove that this method reveals the stable and unstable manifolds of 
hyperbolic points in four particular 2D cases: a hyperbolic saddle point for linear autonomous systems, a hyperbolic saddle point for  nonlinear autonomous systems, a hyperbolic saddle point for linear nonautonomous systems and a hyperbolic saddle point 
for nonlinear nonautonomous systems. We also discuss further rigorous results which show the ability of LDs to highlight additional  invariants sets, such as $n$-tori. These results are just a simple extension of the ergodic partition theory which we illustrate  by applying this methodology to  well-known examples, such as the planar field of the harmonic oscillator and the 3D  ABC flow. Finally, we 
provide a thorough discussion on the requirement  of the objectivity (frame-invariance) property   for tools designed to reveal phase space structures and their implications for Lagrangian descriptors.
$n$-tori, invariant sets.
\end{abstract}

\section{Introduction}
\label{sec:introduction}

Lagrangian descriptors were first introduced in the literature by \cite{chaos} in the form of a function, denoted $M$, that was used to provide a definition for {\em  distinguished trajectories}. The mathematical construction of distinguished trajectories generalized the notion of {\em distinguished hyperbolic trajectory}, first discussed in \cite{kayo}, by including also trajectories with an elliptic type of stability. Distinguished trajectories were highlighted by special minima of the function $M$ referred to as {\em limit coordinates}. 

In the past few years the applicability of the concept of Lagrangian Descriptor  has been extended and has become a method for detecting invariant manifolds of hyperbolic trajectories \citep{prl}. Invariant manifolds were highlighted by  "singular features" of both the function $M$ and some of its generalizations (see \citep{jfm,cnsns}). Since these early papers, numerous applications of Lagrangian descriptors have been given, e.g. in \cite{amism11}, where they were used in the context of atmospheric sciences to reveal the Lagrangian structures that define transport routes across the Antarctic polar vortex. This work was extended in \cite{alvaro2}, where LDs were applied to analyze the Lagrangian structures associated with Rossby wave breaking. In the field of magnetohydrodynamics, Lagrangian descriptors have also been shown to be useful for studying the influence of coherent structures on the saturation of a nonlinear dynamo in \cite{rempel}. There are also several applications in oceanography. In \cite{mmw14}, LDs were used to analyze transport in a region of the Gulf of Mexico relevant to the Deepwater Horizon oil spill. \cite{ggmwm15} have applied this tool to analyze the search strategy for debris from the missing MH370 flight followed by the Australian Maritime Authorities, and recently \cite{ggrmcw16} have studied the role played by LDs in the management of the Oleg Naydenov oil spill that took place in the south of Gran Canaria. All these works are related to various aspects of fluid dynamics. However, Lagrangian descriptors can be applied to the general study of the phase space structure of dynamical systems in different contexts. This has recently been illustrated in several applications of the tool to fundamental problems in chemical reaction dynamics. In particular, it has been applied to a study of  chemical reactions under external time-dependent driving in \cite{CH15}, a study of phase space structure and reaction dynamics for a class of `'barrierless reactions'' in \cite{Jh2015}, and to a study of the isomerization dynamics of ketene in \cite{CH16}.

The ability of LDs to reveal invariant manifolds  has been established  in the references above from a  phenomenological and numerical point of view, however a rigorous framework is missing in these works. Recently \cite{carlos} have provided rigorous proofs in the framework of discrete maps,  where it is  precisely defined what is meant by the phrase "singular features". One of the goals of this article is to  extend those results to continuous time dynamical systems. In order to simplify the demonstrations, this paper provides a new way of constructing Lagrangian descriptors in the same spirit as in \cite{carlos}. The idea is based on considering the $p$-norm of each velocity component, instead of the $p$-norm of the modulus of the velocity. This idea follows the heuristic argument discussed by \cite{cnsns} of integrating positive quantities along particle trajectories, and the positive quantity is  such  that it results in tractable proofs. The choice allows us to mathematically prove that the stable and unstable manifolds of hyperbolic trajectories  in the selected examples are detected as singular features of the Lagrangian descriptor. As  in \cite{carlos}, we are able to make the  notion of ``singular feature'' mathematically precise. 

This paper discusses further rigorous results found in the literature on the ergodic partition theory. These are based on the evaluation of averages along trajectories for obtaining invariant sets (cf. \cite{mezic3, susuki}). We show that LDs are directly related to these findings and thus they also capture coherent structures described  as $n$-tori. To illustrate the full potential of this technique, we apply it to the well-known ABC flow. Finally, we discuss the issue of objectivity and how phase space geometry behaves under coordinate transformations. In this context we show, with some examples from the literature, how outputs of  LDs in different frames consistently reproduce phase portraits.

This paper is organized as follows. In Section \ref{sec:redefinition} a theoretical mathematical background for Lagrangian descriptors is developed, which is based on providing an alternative definition of LDs, and we discuss a variety of Hamiltonian examples that are variations of the linear saddle. Section \ref{sec:non_Ham_sys} presents some results on non-Hamiltonian systems. In Section \ref{sec:LD_elliptic} we illustrate  the link between   Lagrangian descriptors and the ergodic partition theory in the  linear  elliptic case. Section \ref{sec:3D} is devoted to the application of LDs to a well-known 3D 
example, the ABC flow, providing evidence of the effectiveness of LDs in the detection of both invariant manifolds and invariant tori in 3D flows by means of, respectively, singular features and contours of converged averages. Section \ref{sec:LD_obj} offers a detailed discussion of the objectivity  (frame-invariance) property in the context of LDs and the ability of LDs to provide the correct description of phase space structures in different frames, as well as a general consideration of the objectivity property requirement for tools in relation to their capability for revealing Lagrangian structures in phase space. Finally, in Section  \ref{sec:conclusions} we present the conclusions.

\section{Rigorous results for Lagrangian Descriptors}
\label{sec:redefinition}

In this section we provide some rigorous results allowing us to establish a  theoretical framework for Lagrangian descriptors. In order to achieve this goal we propose an alternative definition of LDs for $n$-dimensional vector fields with arbitrary time dependence, following the ideas developed in \cite{carlos} for the discrete time setting. We consider the general time-dependent vector field,
\begin{equation}
  \frac{d\mathbf{x}}{dt} = \mathbf{v}(\mathbf{x},t), \quad \mathbf{x} \in \mathbb{R}^n \;,\; t \in \mathbb{R}
\end{equation}

\noindent where $\mathbf{v}(\mathbf{x},t) \in C^r$ ($r \geq 1$) in $\mathbf{x}$ and continuous in time. The definition of LDs depends on the initial condition $\mathbf{x}_{0} = \mathbf{x}(t_0)$, on the time interval $[t_0-\tau,t_0+\tau]$, and takes the form,
\begin{equation}
M_p(\mathbf{x}_{0},t_0,\tau) = \displaystyle{\int^{t_0+\tau}_{t_0-\tau} \sum_{i=1}^{n} |\dot{x}_{i}(t;\mathbf{x}_{0})|^p \; dt}
\label{M_function}
\end{equation}

\noindent where $p \in (0,1]$ and $\tau \in \mathbb{R}^{+}$ are freely chosen parameters,  and the overdot symbol represents the derivative with respect to time.

\subsection{The autonomous saddle point}
\label{sec:linear_auto_saddle}

The first example that we analyze is the Hamiltonian linear saddle point. The velocity field is given by:
\begin{equation}
\left \{ \begin{array}{cclc}
   \dot{x} & = & \lambda x  &\\
   \dot{y} & = & -\lambda y \\
\end{array}\right . \;\;,\;\; \lambda > 0
\label{eq:lin_aut_saddle_point}
\end{equation}

\noindent For an initial condition $(x_0,y_0)$, the unique solution of this system is:
\begin{equation}
\left \{ \begin{array}{cclc}
   x(t,x_0) & = & x_0 \; e^{\lambda t}  &\\
   y(t,y_0) & = & y_0 \; e^{-\lambda t} \\
\end{array}\right . \;\;,\;\; \lambda > 0
\label{eq:lin_aut_saddle_point_ini_cond}
\end{equation}

\noindent For this example the origin $(0,0)$ is a hyperbolic fixed point with stable and unstable manifolds:
\begin{equation}
  W^s(0,0) = \left\lbrace (x,y) \in \mathbb{R}^2 | x=0, y \neq 0 \right\rbrace,
\end{equation}
\begin{equation}
  W^u(0,0) = \left\lbrace (x,y) \in \mathbb{R}^2 | y=0, x \neq 0 \right\rbrace.
\end{equation}

\noindent For simplicity we assume, without loss of generality, that $t_0 = 0$ (this is possible for autonomous systems) and we apply \eqref{M_function} to \eqref{eq:lin_aut_saddle_point} to obtain:
\begin{equation}
  \begin{array}{ccl}
    M_p((x_0,y_0),t_0,\tau) & = & \displaystyle{\int^{\tau}_{-\tau} |\lambda x_0 e^{\lambda t}|^p + |-\lambda
y_0 e^{-\lambda t}|^p dt } \\
      &   &                \\
      & = & \displaystyle{\left ( |x_0|^p + |y_0|^p \right ) \frac{\lambda^{p-1}(e^{\lambda p \tau}- e^{-\lambda p \tau})}{p}} = 2 \displaystyle{\left ( |x_0|^p + |y_0|^p \right ) \frac{\lambda^{p-1}\sinh (\lambda p \tau)}{p}}. \\
  \end{array}
\label{eq:lagrangian_aut_saddle_point}
\end{equation}

\noindent This expression allows us to conclude the following theorem, which is proven exactly {in the same way as Theorem 1 in \cite{carlos}. 

\begin{theorem}\label{thm:lin_aut_saddle_point}
  Consider a vertical line perpendicular to the unstable manifold of the origin. Then the derivative of $M_p$, $p \leq 1$, along this line does not exist on the unstable manifold of the origin.
\par
  Similarly, consider a horizontal line perpendicular to the stable manifold of the origin. Then the derivative of $M_p$, $p \leq 1$ along this line does not exist on the stable manifold of the origin.
\end{theorem}

This theorem is graphically illustrated in Figure \ref{Saddle_vs_Mp}. Note that this result holds {\bf for any} finite value of $\tau$, i. e., $M_p((x_0,y_{0}),t_0,\tau) < \infty$ possesses, {\it for any} finite $\tau$, singularities along the stable and unstable manifolds. The definition of this singularity is made precise as follows:

\begin{definition}
Given $\tau$ and $p$, an orientable surface $\phi$ with normal vector $\mathbf{n}$ is said to be a singular feature of $M_{p}(\cdot,t_{0},\tau)$, if for every $\mathbf{x}_{0} \in \phi$ the normal derivative $\frac{\partial M_{p}}{\partial \mathbf{n}}(\mathbf{x}_{0},t_{0},\tau)$ does not exist.  \label{def}
\end{definition} 

\begin{remark} We emphasize that the term $| \cdot |^p + |\cdot|^p, \, p \le 1$, at the end of expression (\ref{eq:lagrangian_aut_saddle_point}), whose arguments (denoted by ``$\cdot$'') vanish on the stable and unstable manifolds, is the essential ``structural  feature'' of $M_p$ that gives rise to the ``singular nature'' (i.e. unboundedness or discontinuity in the derivative) of the LD along the stable and unstable manifolds of the hyperbolic trajectory. This feature will appear explicitly in $M_p$ for the benchmark examples to follow.
\end{remark}

\begin{figure}[H]
\centering
{\includegraphics[scale = 0.36]{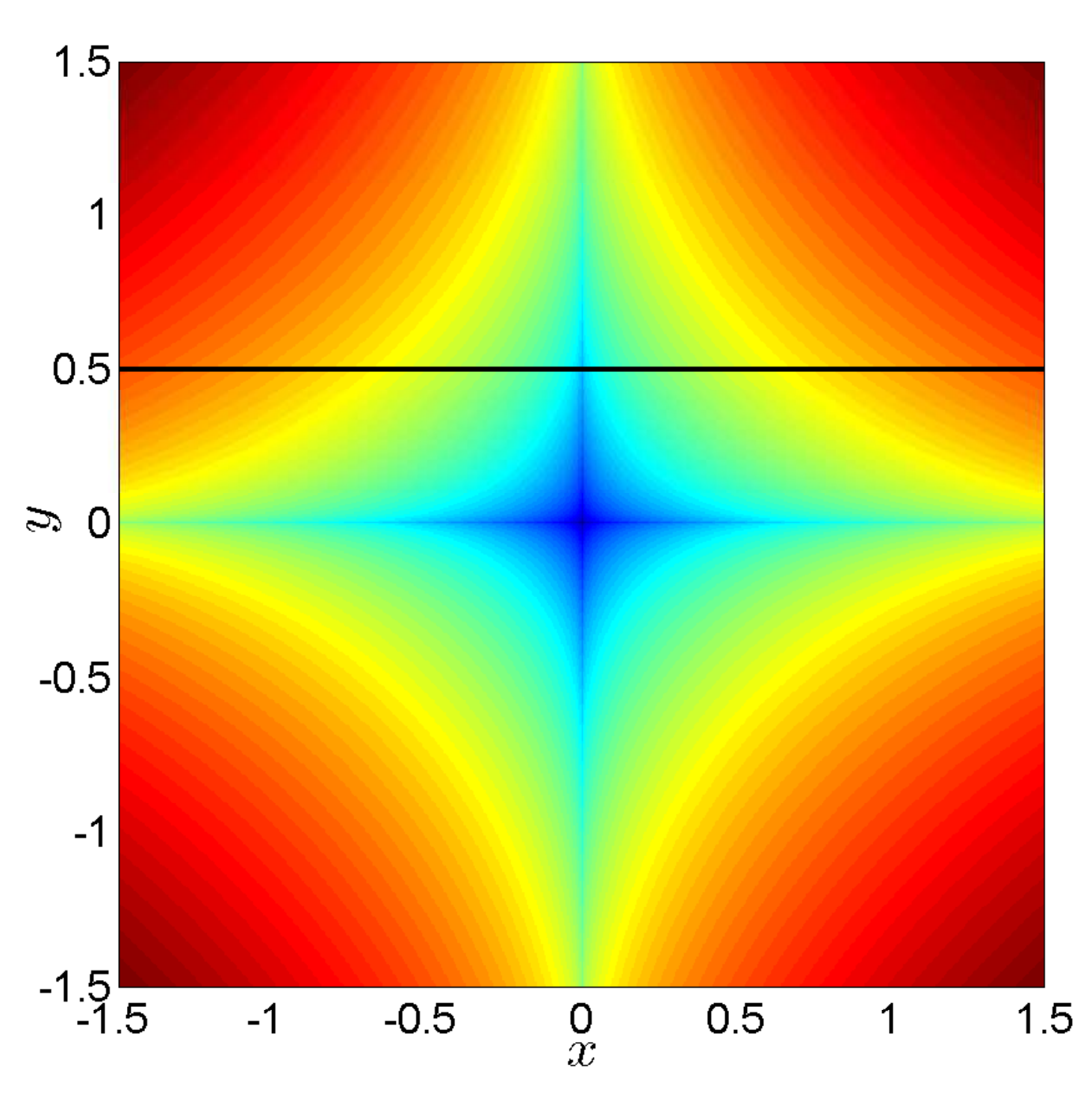}}
{\includegraphics[scale = 0.32]{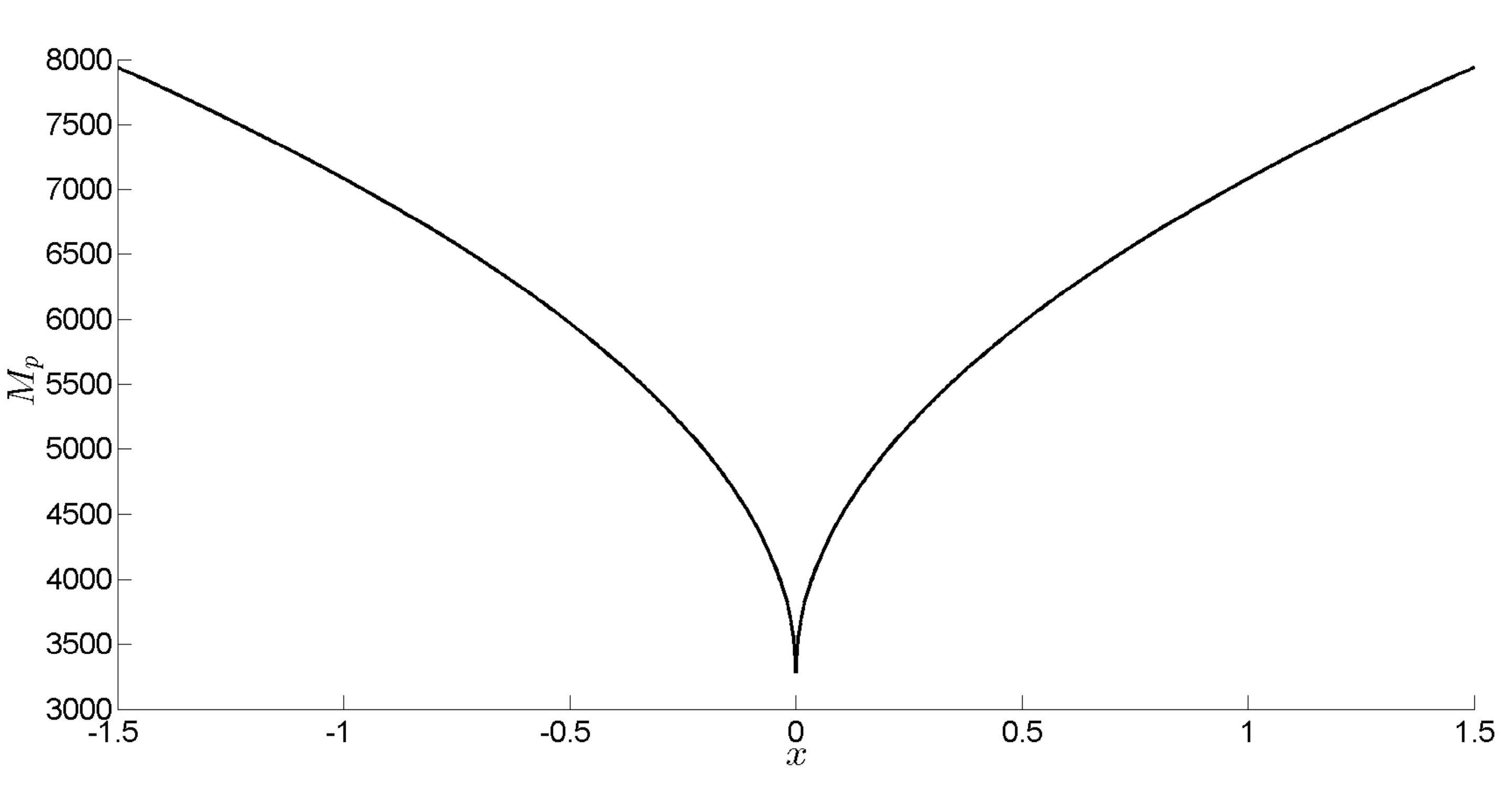}}
\caption{The left-hand panel shows contours of $M_{p=0.5}$ for system (\ref{eq:lin_aut_saddle_point}) with $\tau=15$ and $\lambda=1$. The contours of $M_p$ are computed on a grid with spacing $0.005$ and the integration time step of the vector field is chosen to be $0.1$. The horizontal black line is at $y=0.5$. The right-hand panel shows the graph of $M_p$ along this line, which illustrates the singular nature of the derivative of $M_p$ on the stable manifold $\lbrace x=0 \rbrace$.}
\label{Saddle_vs_Mp}
\end{figure}

\subsection{The rotated saddle point}
\label{sec:rotated_saddle}

In the previous example Lagrangian descriptors were shown to be singular along the stable and unstable manifolds of the hyperbolic fixed point at the origin for any finite value of $\tau$. However for the examples given in \cite{prl} and \cite{cnsns}, it was shown that a sufficiently large $\tau$ was required in order to obtain ``sharp'' images/figures  on/over the manifolds. In the next example we illustrate this particular role of $\tau$ by considering the same example of the previous section, but rotated 45$\degree$. The autonomous system corresponding to the 45$\degree$ rotated saddle has the form:
\begin{equation}
\left \{ \begin{array}{ccl}
   \dot{x} & = & y \\
   \dot{y} & = & x \\
\end{array}\right . \;\;,
\label{eq:rotated_saddle_point}
\end{equation}

\noindent The solution of this system yields,
\begin{equation}
\left \{ \begin{array}{ccl}
   x(t) & = & x_0\cosh t + y_0\sinh t  \\
   y(t) & = & x_0\sinh t + y_0\cosh t  \\
\end{array}\right . \;\;,
\label{eq:hyperbolic_solution}
\end{equation}

\noindent or equivalently,
\begin{equation}
\left \{ \begin{array}{ccl}
   x(t) & = & ae^t + be^{-t} \\
   y(t) & = & ae^t - be^{-t} \\
\end{array}\right . \;\;,
\label{eq:exponential_solution}
\end{equation}

\noindent where 
\begin{equation}
  a = \frac{x_0+y_0}{2} \;\;,\;\; b = \frac{x_0-y_0}{2} \;.
\end{equation}
\noindent
Note that $a=0$ corresponds to the stable manifold and $b=0$ corresponds to the unstable manifold.

The function $M_p$ for this example takes the form:
\begin{equation}
  M_p((x_0,y_0),t_0,\tau) = \displaystyle{\int^{\tau}_{-\tau} |ae^t - be^{-t}|^p + |ae^t + be^{-t}|^p dt}. 
\end{equation}

\noindent In this example, it is not possible to compute analytically the integrals which define $M_p$. However, we are able to compute approximations of $M_p$ that are sufficiently accurate to enable the understanding of the relationship between singularities of $M_p$ and the stable and unstable manifolds for both small and large $\tau$ limits.

First we study the behavior of $M_p$ in order to show where the singularities of the derivative of $M_p$ appear for small $\tau=\tau_0$. The Taylor expansions of $\sinh t$ and $\cosh t$ are:
\begin{equation}
  \sinh t = t + \frac{t^3}{3!} + \frac{t^5}{5!} + \cdots \;\;,\;\; \cosh t = 1 + \frac{t^2}{2!} + \frac{t^4}{4!} + \cdots
\end{equation}

\noindent Using \eqref{eq:hyperbolic_solution} and the Taylor expansion of $M_p((x_0,y_0),t_0,\tau_0)$ at $\tau_0=0$,  

\begin{equation}
\begin{array}{ccl}
  M_p((x_0,y_0),t_0,\tau_0) & = & M_p((x_0,y_0),t_0,0) + \displaystyle{\frac{\partial M_p((x_0,y_0),t_0,0)}{\partial \tau_0}}\tau_0 
+ O(\tau_0^2) \\
                      &   & \\
                      & = & 0 + \left( |x_0\sinh \tau_0 + y_0\cosh \tau_0|^p + |x_0\cosh \tau_0 + y_0\sinh 
\tau_0|^p \right )_{\tau_0=0}\tau_0  \\
                      &   & \\
                      & + & \left( |x_0\sinh (-\tau_0) + y_0\cosh (-\tau_0)|^p + |x_0\cosh (-\tau_0) + y_0\sinh 
(-\tau_0)|^p \right )_{\tau_0=0}\tau_0 + O(\tau_0^2)\\
		      &   & \\
                      & \sim & 2 \left( |y_0|^p + |x_0|^p \right) \tau_0
\end{array}
\end{equation}

\noindent
where in the last step we have used that  $\sinh \tau_0 =-\sinh (-\tau_0) =0$, and that  $\cosh \tau_0 =\cosh (-\tau_0) =1$. Therefore the singularities of the derivative of 
 $M_p$ appear on the lines $x=0$ and $y=0$ for small $\tau$, which are not the stable and unstable manifolds of the origin. 

Now we consider the case of large values of $\tau$, and fixed $\tau_0$. In order to analyze this case we divide the integral into three parts,
\begin{equation}\label{eq:M_expansion}
M_p((x_0,y_0),t_0,\tau) = \displaystyle{M_p((x_0,y_0),t_0,\tau_0) + \int^{\tau}_{\tau_0} |\dot{x}(t)|^p + |\dot{y}(t)|^p dt + 
\int^{-\tau_0}_{-\tau} |\dot{x}(t)|^p + |\dot{y}(t)|^p dt}.
\end{equation}

\noindent Since $\tau_0$ is fixed and $\tau$ is large enough, we can expand the last part of \eqref{eq:M_expansion} to yield, 
\begin{equation}\label{eq:M_assymptotic_expansion_plus}
\begin{array}{rcl}
  \displaystyle{\int^{\tau}_{\tau_0} |\dot{x}(t)|^p + |\dot{y}(t)|^p dt} & = & \displaystyle{\int^{\tau}_{\tau_0}\left 
| ae^t - be^{-t}\right |^p + \left |ae^t + be^{-t}\right |^p dt }\\
                                                                         &   &  \\
                                                                         & = &  
2 \displaystyle{\int^{\tau}_{\tau_0} |ae^t|^p + O({|b|}/{|a|^{1-p} e^{2 t -p }}) \; dt. }
\end{array}
\end{equation}

\noindent Analogously, for the range of negative values we obtain,  
\begin{equation}\label{eq:M_assymptotic_expansion_minus}
\begin{array}{rcl}
  \displaystyle{\int^{-\tau_0}_{-\tau} |\dot{x}(t)|^p + |\dot{y}(t)|^p dt} & = & 
\displaystyle{\int^{-\tau_0}_{-\tau}\left 
| ae^t - be^{-t}\right |^p + \left |ae^t + be^{-t}\right |^p dt }\\
                                                                           &   &  \\
                                                                           & = &  
2 \displaystyle{\int^{-\tau_0}_{-\tau} |be^{-t}|^p + O(|a|/|b|^{1-p}  e^{-2t+p}) \; dt. } 
\end{array}
\end{equation}

\noindent {Using the fact that $\tau\gg\tau_0> 0$, it is clear that the leading order terms of \eqref{eq:M_assymptotic_expansion_plus} and \eqref{eq:M_assymptotic_expansion_minus} after the integration are given by
\begin{equation}
  2    \frac{|a|^p}{p}  e^{\tau p}   \qquad \text{and} \qquad    2    \frac{|b|^p}{p}  e^{\tau p}  
\end{equation}}
\noindent
respectively. Consequently, 
\begin{eqnarray}
  M((x_0,y_0),t_0,\tau) & = & M((x_0,y_0),t_0,\tau_0) + \displaystyle{\int^{\tau}_{\tau_0}\left | ae^t - be^{-t}\right |^p + 
\left |ae^t + be^{-t}\right |^p dt }\nonumber\\
                    &   & \nonumber  \\
                    & + & \displaystyle{\int^{-\tau_0}_{-\tau}\left | ae^t - be^{-t}\right |^p + 
\left |ae^t + be^{-t}\right |^p dt }\\
		                        &   & \nonumber  \\
                    & = & 2 (|x_0|^p + |y_0|^p) \tau_0 + \displaystyle{\frac{2^{1-p} e^{\tau p}}{p} 
\cdot \left ( |x_0+y_0|^p + |x_0-y_0|^p\right )+B.} \label{eq21}
\end{eqnarray}

\noindent In (\ref{eq21}) $B$ depends on $\tau_0$ and lower order terms in $\tau$. The singularities previously discussed, observed for small enough $\tau$ along the horizontal and vertical axis, are still present in the first term, although the weight of this term makes it negligible when compared to the second term. For large $\tau$, (\ref{eq21}) depicts the singular features of $M_p$ over/at the lines $\lbrace y=-x \rbrace$
(the stable manifold) and $\lbrace y=x \rbrace$ (the unstable manifold). The evolution of the contours of $M_p$ and its singularities can be seen in Figure \ref{M_p_sequence} and \ref{M_p_singularities}, respectively, as we increase the  integration time parameter
$\tau$. It is clear from these figures that the longer we integrate the system (increasing $\tau$) the closer we get to patterns that enhance the diagonal structure.
For very small $\tau$, the major contribution of $M_p$ comes from the first term in (\ref{eq21}), while for intermediate $\tau$ values the contribution of $B$ is the dominant one, although we do not know its explicit expression. Finally, for very large $\tau$, the dominant term in (\ref{eq21}) is the second one. We remark that Figures \ref{M_p_sequence} and \ref{M_p_singularities} show that the very large $\tau$  required
for the singular features of $M_p$ to be visibly aligned along the  stable and unstable manifolds, is achieved already at $\tau=5$, {\em i.e.} a {\em finite} $\tau$ value at which of course $M_p<\infty$. Eq. (\ref{eq21}) confirms the ability of $M_{p}$ to highlight invariant manifolds aligned in different directions by means of singular features as defined in \ref{sec:linear_auto_saddle}.

\begin{figure*}[htbp!]
  \centering
  \subfigure[$M_p$ for $\tau=0.005$]{\includegraphics[width=0.42\linewidth]{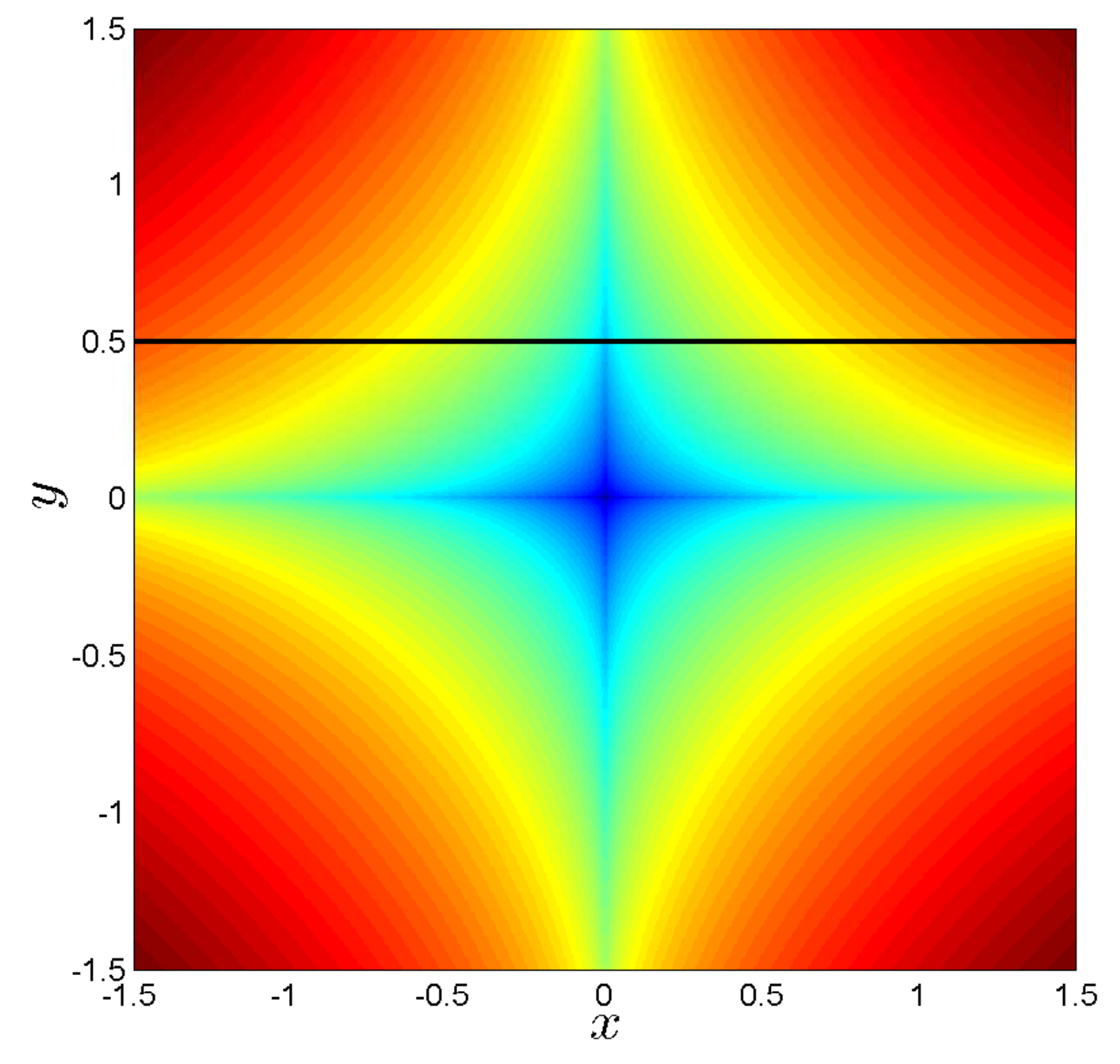}}
  \subfigure[$M_p$ for $\tau=1$]{\includegraphics[width=0.42\linewidth]{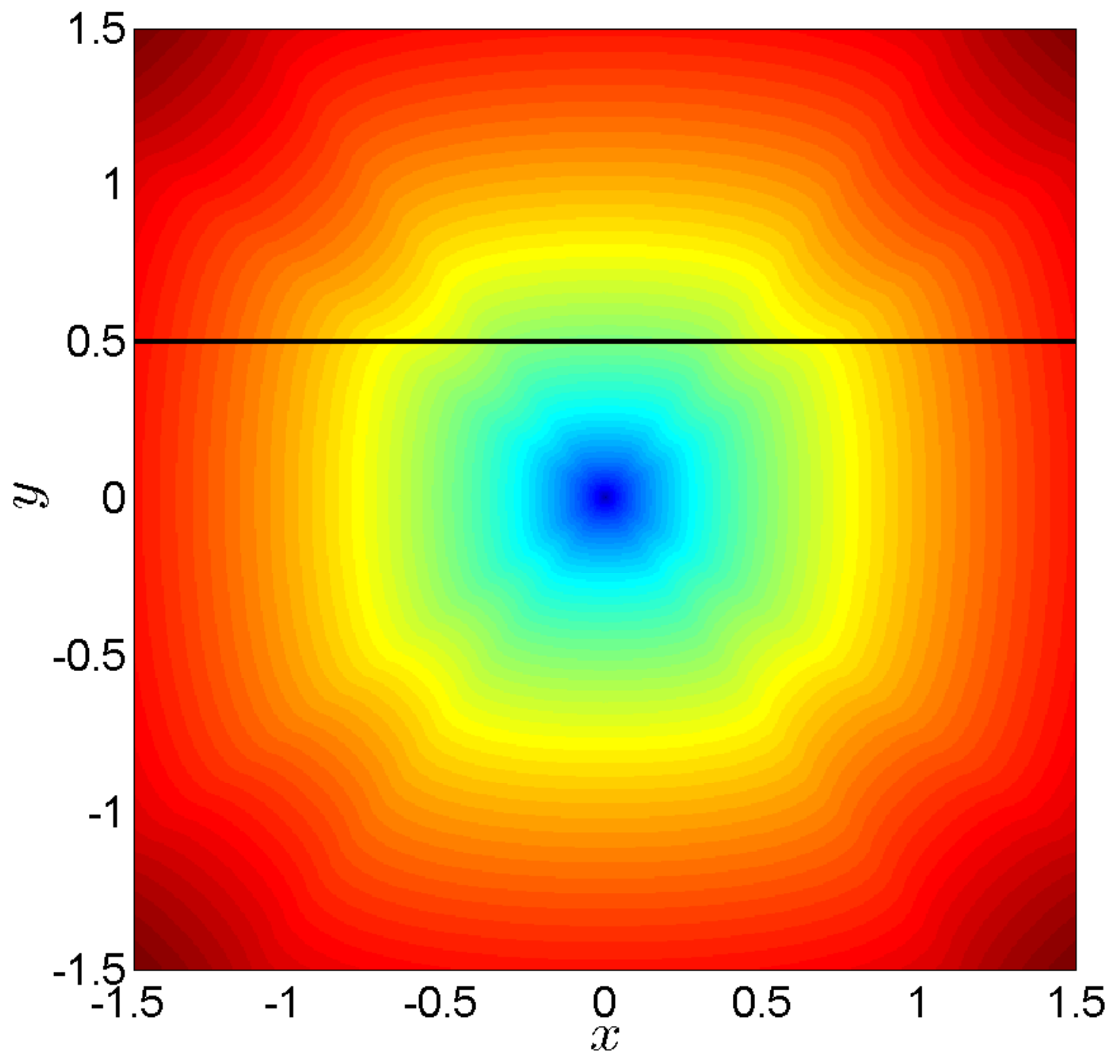}} \\
  \subfigure[$M_p$ for $\tau=2.5$]{\includegraphics[width=0.42\linewidth]{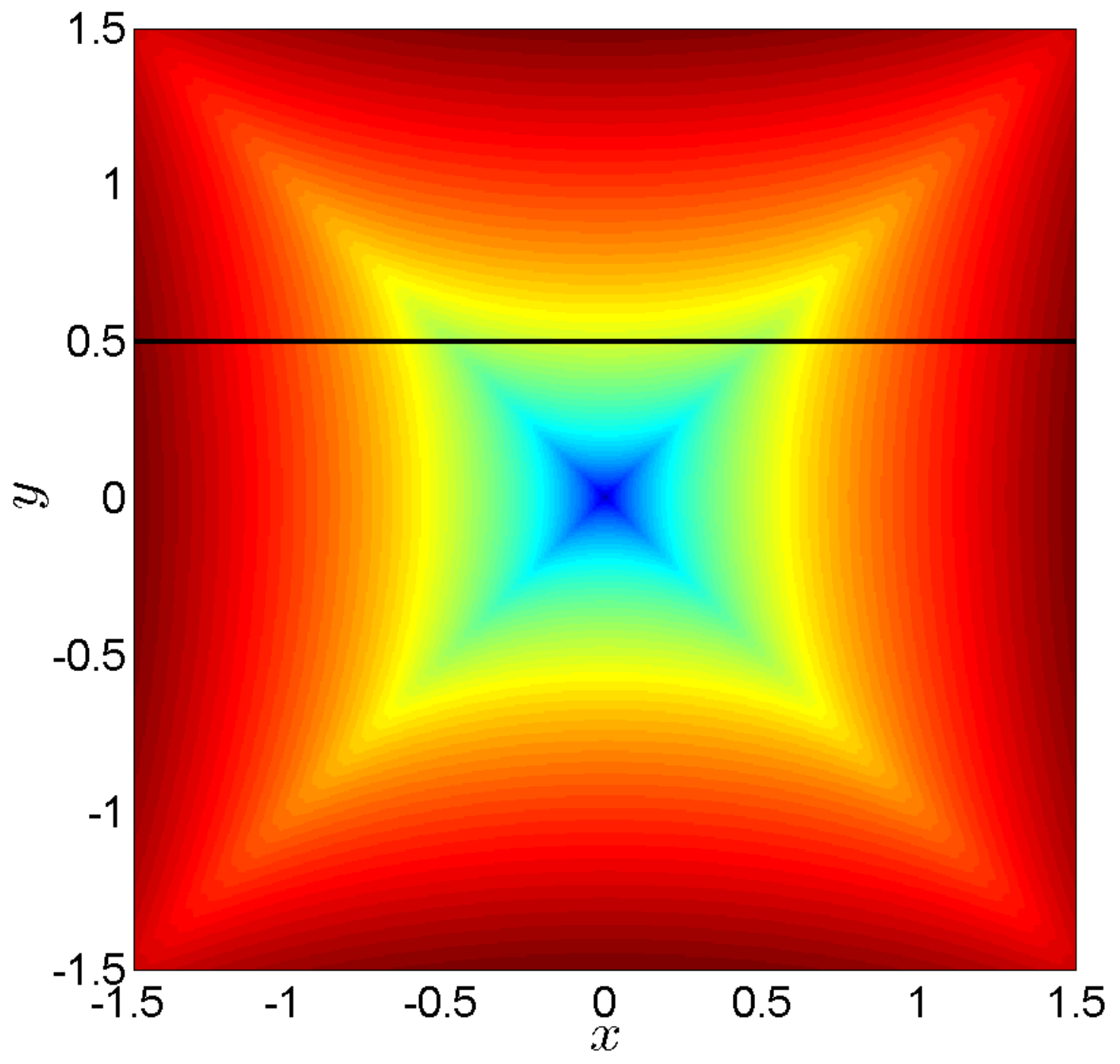}}
  \subfigure[$M_p$ for $\tau=5$]{\includegraphics[width=0.42\linewidth]{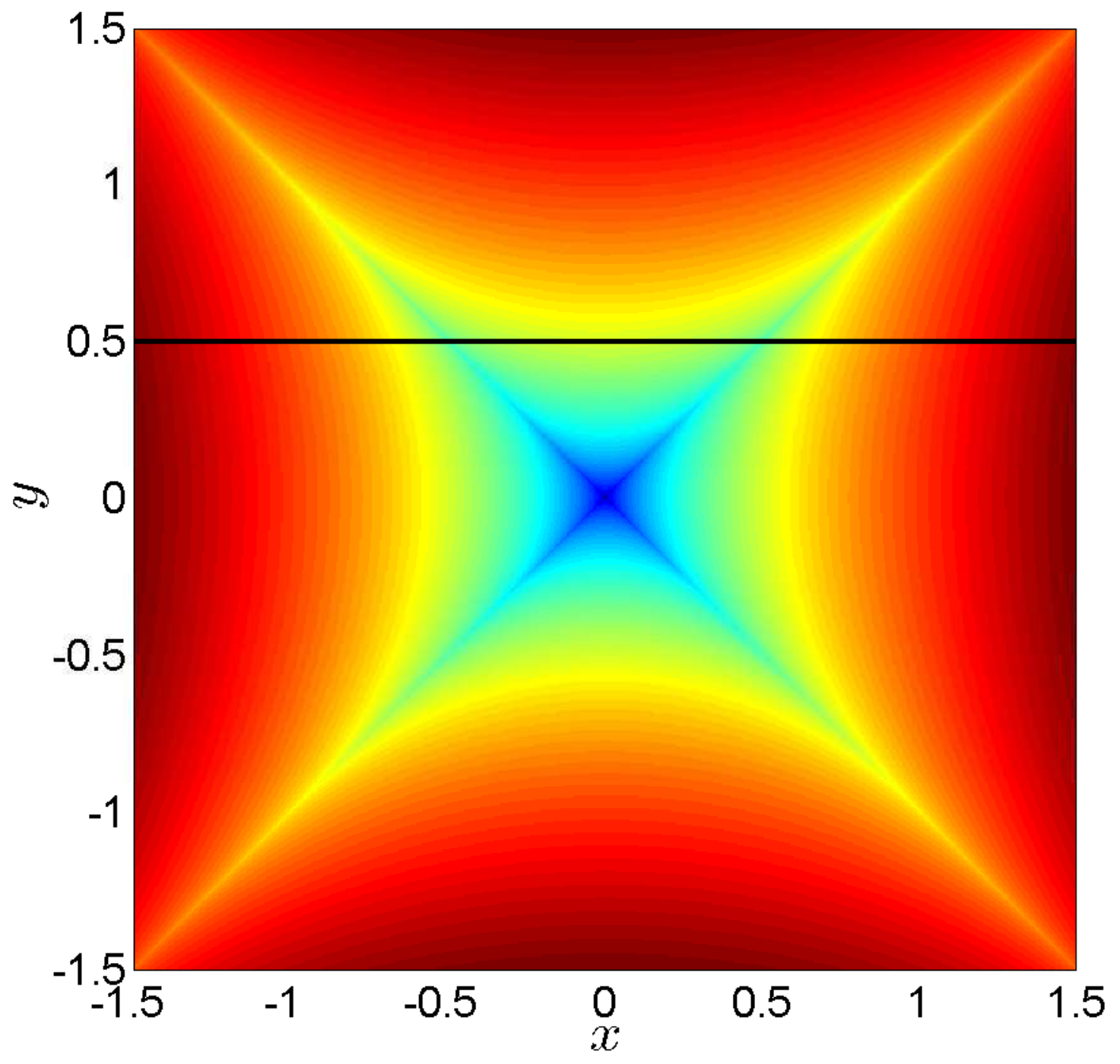}} \\
  \caption{$M_p$ function for $p = 0.5$ and using different values of $\tau$. For this example the integration time step 
and grid spacing are 0.005. The black line is defined at $y=0.5$ and is used in Figure \ref{M_p_singularities}.}
  \label{M_p_sequence}
\end{figure*}

\begin{figure*}[htbp!]
  \centering
  \subfigure[$M_p$ for $\tau=0.005$]{\includegraphics[width=0.49\linewidth]{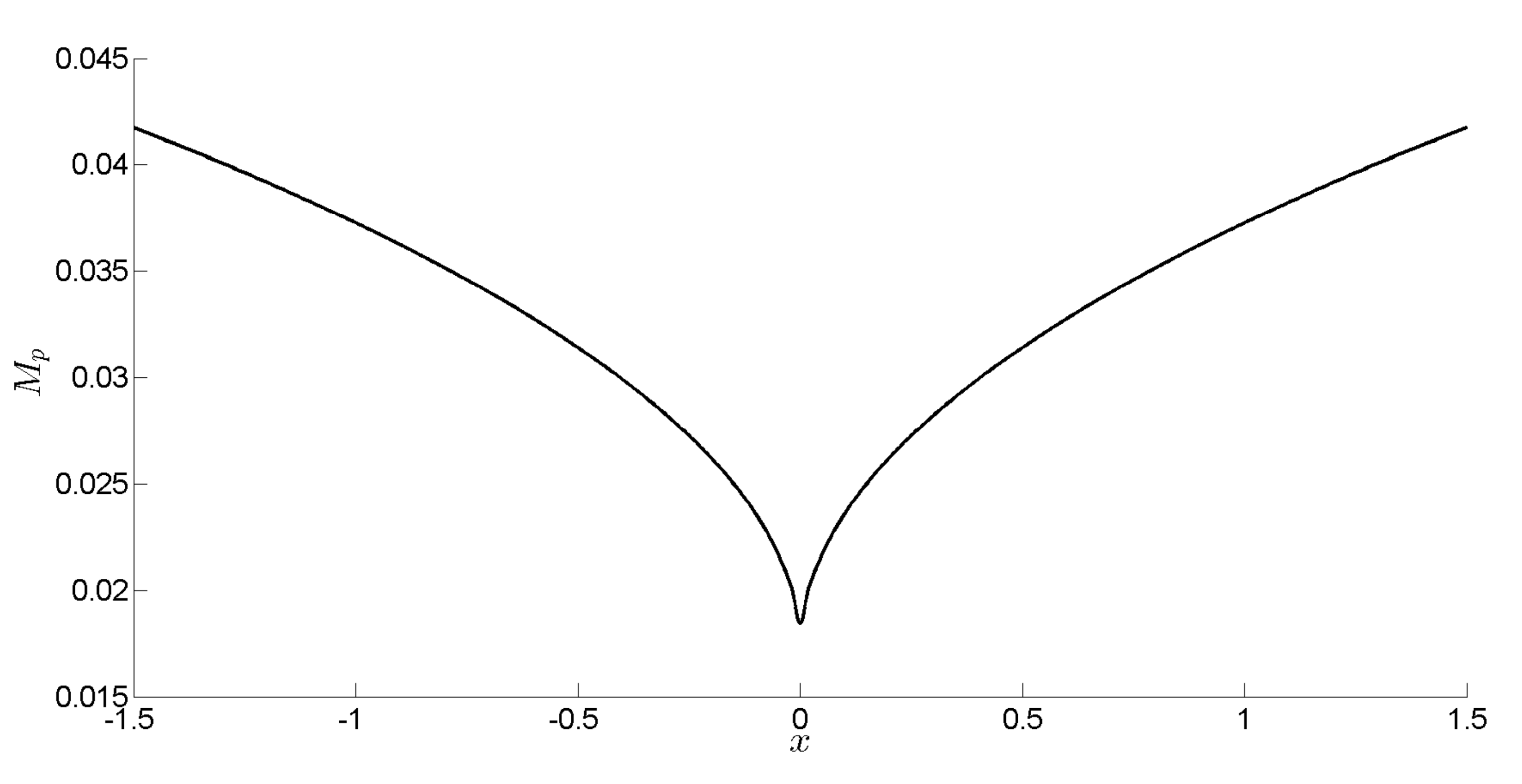}}
  \subfigure[$M_p$ for $\tau=1$]{\includegraphics[width=0.49\linewidth]{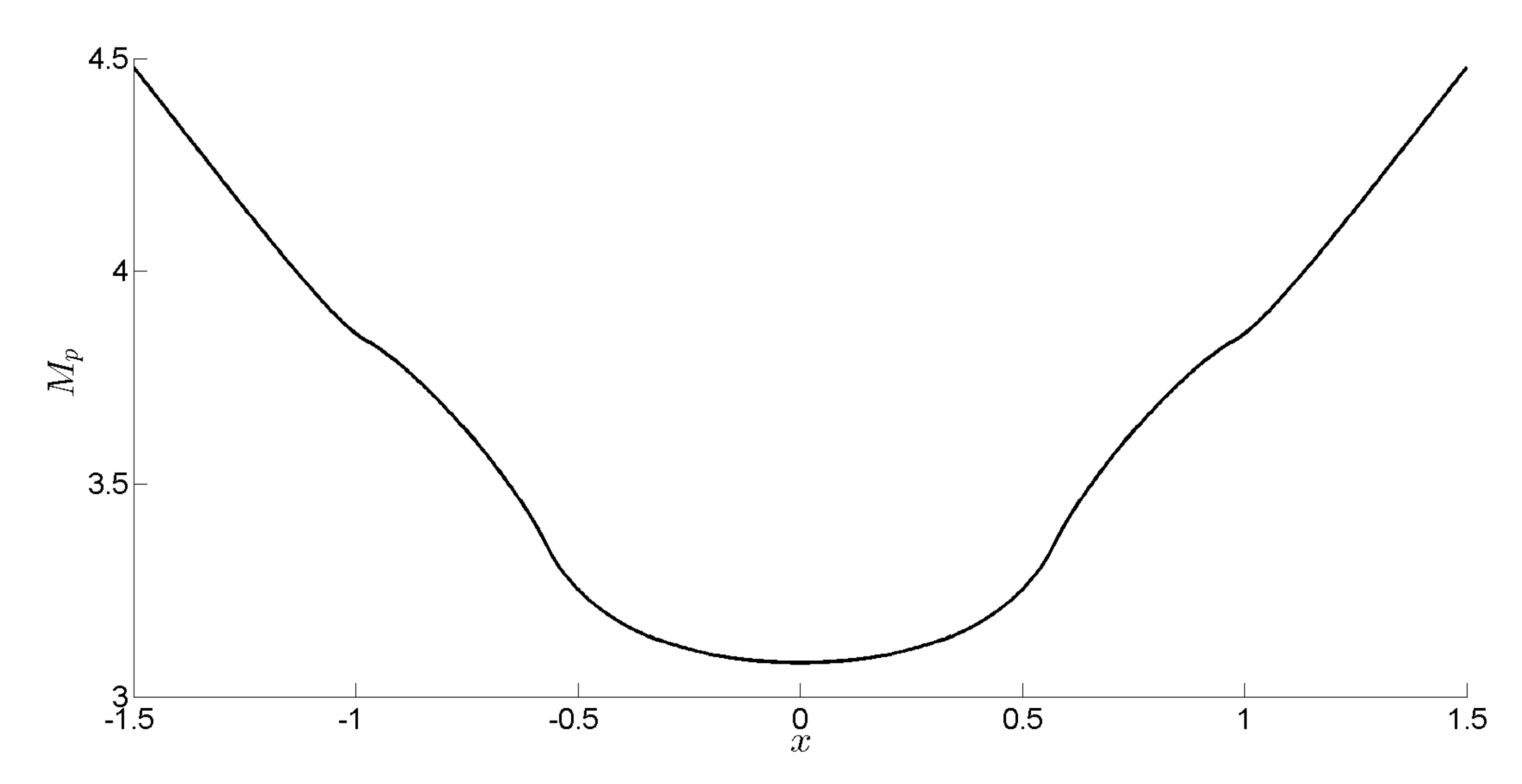}}\\
  \subfigure[$M_p$ for $\tau=2.5$]{\includegraphics[width=0.49\linewidth]{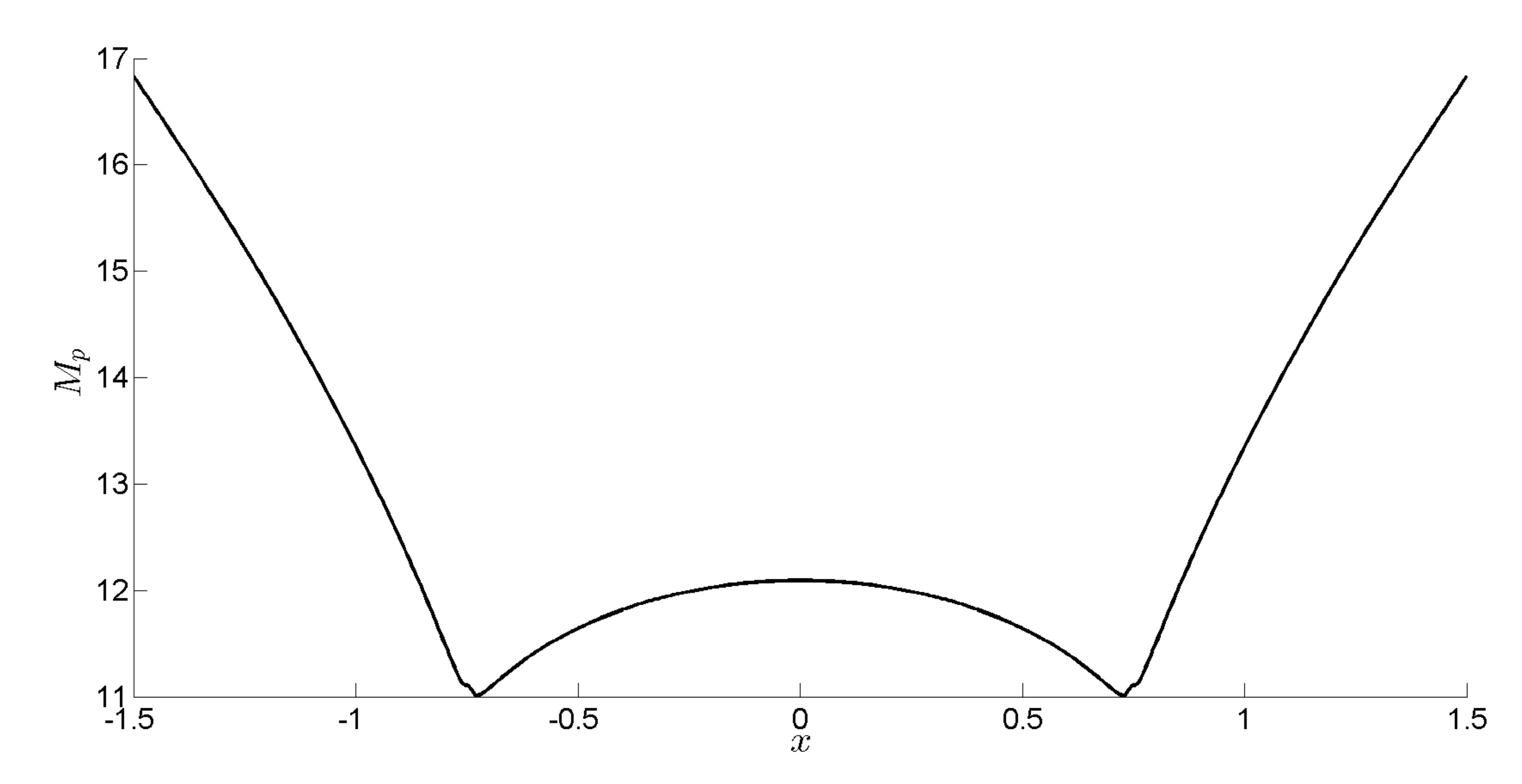}}
  \subfigure[$M_p$ for $\tau=5$]{\includegraphics[width=0.49\linewidth]{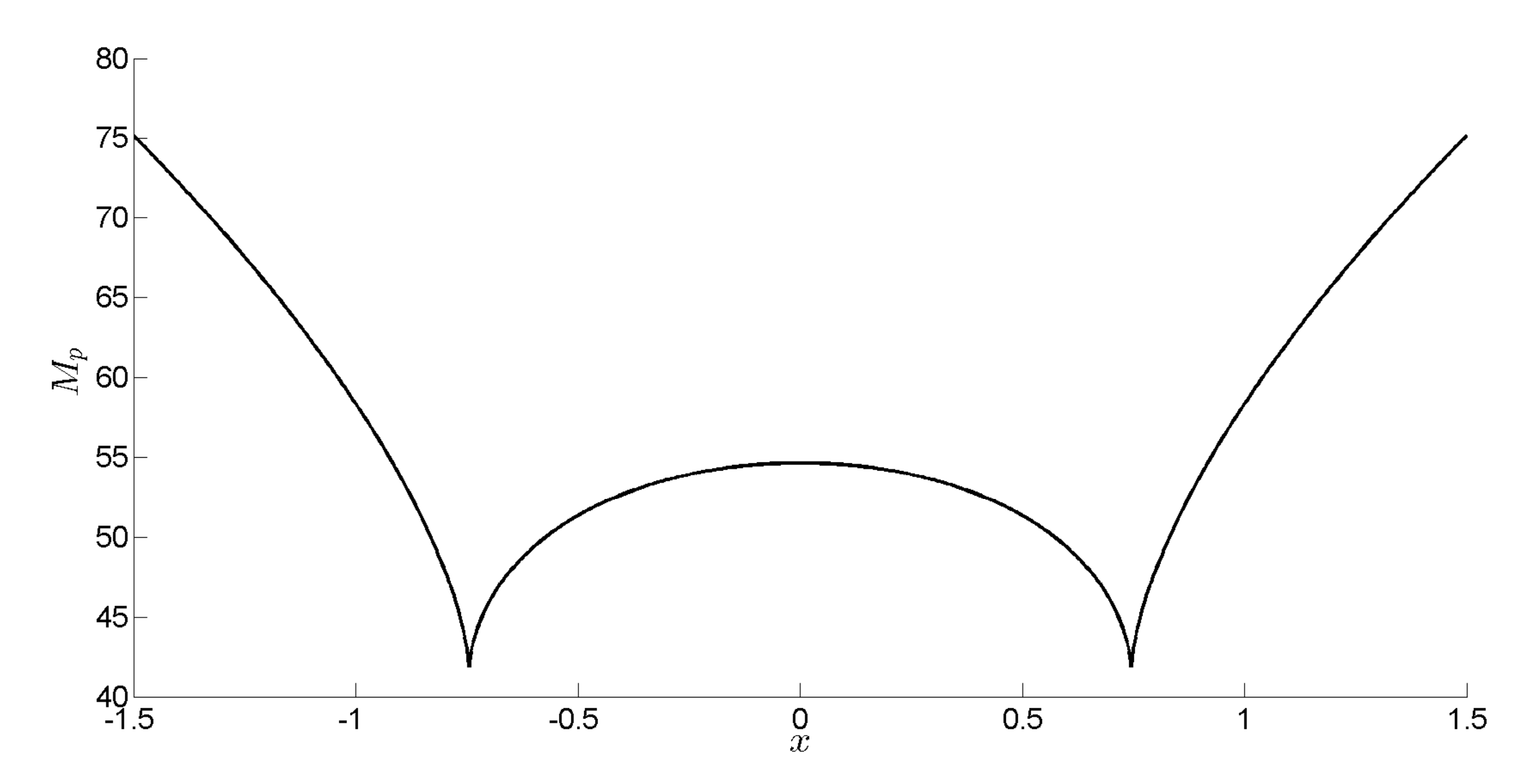}}\\
  \caption{Evolution of the singular features of $M_p$ for $p = 0.5$ along the line $\lbrace y=0.5 \rbrace$  using different values of $\tau$. }
  \label{M_p_singularities}
\end{figure*}

\subsection{The autonomous nonlinear saddle point}
\label{sec:nonlinear_auto_saddle}

In this section we treat the autonomous nonlinear saddle point by using a theorem by \cite{moser56}. Moser's theorem applies to analytic two-dimensional symplectic maps having a hyperbolic fixed point or, similarly, to two-dimensional time-periodic Hamiltonian vector fields having a hyperbolic periodic orbit (which can be reduced to the former case considering a Poincar\'e map). In this case we are considering the flow of an autonomous Hamiltonian system, which is a one-parameter family of symplectic maps, and therefore Moser's theorem applies.

We consider a two-dimensional autonomous analytic Hamiltonian vector field of the form,
\begin{equation}
\left \{ \begin{array}{cclc}
   \dot{x} & = & H_y(x,y), \\
   \dot{y} & = & -H_x(x,y), \\
\end{array}\right .
\label{eq:nonlin_aut_saddle_point}
\end{equation}

\noindent having a hyperbolic fixed point at the origin. Moser's theorem proves the existence of an area-preserving change of variables,
\begin{equation}
\left \{ \begin{array}{ccl}
   x & = & x(\xi, \eta) \\
   y & = & y(\xi, \eta) \\
\end{array}\right .
\label{eq:moser_change_variables}
\end{equation}

\noindent with inverse,
\begin{equation}
\left \{ \begin{array}{ccl}
   \xi  & = & \xi(x, y) \\
   \eta & = & \eta(x, y) \\
\end{array}\right .
\label{eq:moser_inverse_change_variables}
\end{equation}

\noindent by which the system \eqref{eq:nonlin_aut_saddle_point} is transformed into the following normal form,
\begin{equation}
\left \{ \begin{array}{ccl}
   \dot{\xi}  & = & F_{\eta} \\
   \dot{\eta} & = & -F_{\xi} \\
\end{array}\right .
\label{eq:moser_canonical_change_variables}
\end{equation}

\noindent where $F = F(\xi \eta) = a_0 \xi \eta + a_1 (\xi \eta)^2 + \cdots$ depends only on the product $\xi \eta$ and $a_0 \equiv
\lambda \in \mathbb{R}$, $\lambda \neq 0$  (onwards it will be taken that $\lambda > 0$). It is straightforward to verify that $\frac{d}{dt} (\xi \eta) =0$, or equivalently, $\xi_0 \eta_0 = \xi \eta$. Moreover, if we define $\frac{dF}{dz} (z) = F' (z)$, with $z \equiv \xi \eta$ then \eqref{eq:moser_canonical_change_variables} takes the form,
\begin{equation}
\left \{ \begin{array}{ccl}
   \dot{\xi}  & = & F^{'}\xi \\
   \dot{\eta} & = & -F^{'}\eta \\
\end{array}\right .
\end{equation}

\noindent where $F^{'} = F^{'}(\xi_0\eta_0)$ is constant on trajectories. This last system is easily integrated and its solutions are given by the expressions,
\begin{equation}
\left \{ \begin{array}{ccl}
   \xi  & = & \xi_0e^{F^{'}t}\\
   \eta & = & \eta_0e^{-F^{'}t} .\\
\end{array}\right .
\end{equation}

\noindent Applying the descriptor $M_p$ to this system (and setting $t_0=0$ since it is autonomous) we obtain,
\begin{equation}
\begin{array}{ccl}
  M_p((\xi_0, \eta_0), t_0, \tau)  & = & \displaystyle{\int^{\tau}_{-\tau} |\xi_0e^{F^{'}t}F^{'}|^p +
|\eta_0e^{-F^{'}t}F^{'}|^p dt } = 2 \left ( |\xi_0|^p + |\eta_0|^p\right ) \displaystyle{\frac{|F^{'}|^p}{pF^{'}}}
\sinh (pF^{'}\tau). \\
\end{array}
\label{LDNLHamaut}
\end{equation}

\noindent As we can see above, the derivative of $M_p$ has singularities through the stable $\lbrace \xi = 0 \rbrace$ and unstable $\lbrace \eta =
0 \rbrace$ manifolds in the $\xi$ - $\eta$ coordinates, i.e. the directional derivatives of $M_p$ across $\xi=0$ and $\eta=0$ do not exist on the manifolds.

In order to complete this example, there are two technical points that we must address. From \eqref{LDNLHamaut} we observe that there might be possible singularities in the case when $F^{'}$ vanishes. However, recall that $F'$ is,
\begin{equation}
   \frac{d F}{d z} = \lambda + 2a_1 z + 3a_2 z^2 \cdots,
\end{equation}
  
\noindent where $z \equiv \xi \eta$. For a sufficiently small neighborhood of the origin (no larger than the domain of validity of the normal form) the term $\lambda$ is the dominant term in the series, and therefore $F'$ is not zero.

We have shown that $M_p$ is singular on the stable and unstable manifolds in the $\xi-\eta$ coordinates. However, the  system \eqref{eq:nonlin_aut_saddle_point} was originally expressed in the $x-y$ coordinates. Now we show that $M_p$ is singular on the stable and unstable manifolds in these coordinates. We will carry out the proof for the stable manifold. The proof for the unstable manifold is completely analogous.

For this purpose we use  expression \eqref{eq:moser_change_variables}  and that the stable manifold is given by $\xi=0$ to obtain the following parametrization for the stable manifold in the $x-y$ coordinates,
\begin{equation}
\label{eq:parametrization}
\eta \in \mathbb{R} \quad \longmapsto \quad (x(0,\eta ),y(0,\eta )) \in \mathbb{R}^{2},
\end{equation}

\noindent where $ \left( \frac{\partial x}{\partial \eta}(0,\eta ) , \frac{\partial y}{\partial \eta}(0,\eta ) \right) $ is the tangent vector of this curve at every point $(x(0,\eta ),y(0,\eta ))$. Additionally, since the Jacobian of the transformation (\ref{eq:moser_change_variables}) is non-zero,
\begin{equation}
J(\xi ,\eta )= \frac{\partial x}{\partial \xi} \frac{\partial y}{\partial \eta} (\xi ,\eta ) - \frac{\partial 
y}{\partial \xi} \frac{\partial x}{\partial \eta} (\xi ,\eta ) \not = 0 \quad \text{for all} \quad (\xi ,\eta ) 
\in \mathbb{R}^{2},
\end{equation}

\noindent the pair of vectors
\begin{equation}
\left \{ \left( \frac{\partial x}{\partial \xi}(0,\eta ) , \frac{\partial 
y}{\partial \xi}(0,\eta ) \right) , \left( \frac{\partial x}{\partial \eta}(0,\eta ) , \frac{\partial y}{\partial 
\eta}(0,\eta ) \right) \right \}_{\eta \in \mathbb{R}} \label{vect}
\end{equation}

\noindent cannot be parallel, and thus they form a basis of $\mathbb{R}^{2}$ at every point $(x(0,\eta ),y(0,\eta ))$ of the stable manifold. Since the change of variables (\ref{eq:moser_change_variables}) is analytic, one can construct a family of unit normal vectors $\lbrace \vec{n}(\eta ) \rbrace_{\eta \in \mathbb{R}}$ to the stable manifold (\ref{eq:parametrization}) which can be expressed as follows,
\begin{equation}
\vec{n}(\eta ) = a(\eta ) \left( \begin{array}{c} \frac{\partial x}{\partial \xi}(0,\eta ) \\ 
\frac{\partial y}{\partial \xi}(0,\eta ) \end{array} \right) + b(\eta ) \left( \begin{array}{c} \frac{\partial 
x}{\partial \eta}(0,\eta ) \\ \frac{\partial y}{\partial \eta}(0,\eta ) \end{array} \right)
\end{equation}

\noindent where  $a,b \in C^{0}(\mathbb{R})$ are two scalar functions. Here $a(\eta ) \not = 0, \forall \eta \in \mathbb{R}$, since,
by definition, every normal vector is perpendicular to the tangent vector at every single point of a $C^{1}$-curve, and thus must have a non-zero component along a vector which is not parallel to the tangent direction. On the other hand,  the Jacobian of the transformation only confirms that the vectors (\ref{vect}) are not parallel, but not that they are perpendicular. For this reason $b(\eta )$ can take values different from zero.

Now the directional derivative of $M_{p}$ in the direction $\vec{n}(\eta )$ is given by,
$$\left( \frac{\partial M_{p}}{\partial x} (x(0,\eta ),y(0,\eta )) , \frac{\partial M_{p}}{\partial y} (x(0,\eta 
),y(0,\eta )) \right) \cdot \vec{n}(\eta ) $$
$$= a(\eta ) \left[ \frac{\partial M_{p}}{\partial x} (x(0,\eta ),y(0,\eta )) \frac{\partial x}{\partial 
\xi}(0,\eta ) + \frac{\partial M_{p}}{\partial y} (x(0,\eta ),y(0,\eta )) \frac{\partial y}{\partial \xi}(0,\eta ) 
\right] $$
$$+ b(\eta ) \left[ \frac{\partial M_{p}}{\partial x} (x(0,\eta ),y(0,\eta )) \frac{\partial x}{\partial 
\eta}(0,\eta ) + \frac{\partial M_{p}}{\partial y} (x(0,\eta ),y(0,\eta )) \frac{\partial y}{\partial \eta}(0,\eta 
) \right]$$
\begin{equation}
= a(\eta ) \frac{\partial M_{p}}{\partial \xi}(0,\eta ) + b(\eta ) \frac{\partial M_{p}}{\partial 
\eta}(0,\eta ), \quad \forall \eta \in \mathbb{R}.
\end{equation}

\noindent Since $a(\eta )$ cannot be equal to zero, this derivative is unbounded in the direction of the normal vector  
$\vec{n}(\eta)$ for $p \in (0,1)$ or has a discontinuity for $p=1$. In both cases we are speaking in terms of the
non-differentiability of function $M_p$, therefore displaying a singular feature over the manifolds expressed in the $x-y$ coordinates. Indeed its directional derivative is bounded and continuous only when evaluated with respect to the tangent vector 
direction.

\subsection{The nonautonomous linear saddle point}
\label{sec:linear_nonauto_saddle}

In this section we consider the nonautonomous linear saddle point. Given a 
function $f \in C^1(t_0-\tau, t_0+\tau)$, and $f(t) > 0$ for every $t \in [t_0-\tau, t_0+\tau]$, we define the following vector field,
\begin{equation}
\left \{ \begin{array}{ccl}
   \dot{x} & = & f(t) x,  \\
   \dot{y} & = & -f(t) y, \\
\end{array}\right .
\label{eq:lin_non_aut_saddle_point}
\quad 
\end{equation}

\noindent For any initial condition $(x(t_0),y(t_0))=(x_0,y_0)$, the solution of  \eqref{eq:lin_non_aut_saddle_point} is given by,
\begin{equation}
\left \{ \begin{array}{ccl}
   x(t,x_0) & = & x_0e^{F(t)}  \\
            &   &              \\
   y(t,y_0) & = & y_0e^{-F(t)} \\
\end{array}\right .
\label{eq:lin_non_aut_saddle_point_ini_cond}
\end{equation}

\noindent where $F(t) = \displaystyle{\int^{t}_{0} f(s)ds}$. This system has a stationary hyperbolic trajectory at the origin, with stable manifold given by $W^{s}(0,0)=\lbrace (x,y) : x=0, y \not = 0 \rbrace$ and unstable manifold given by $W^{u}(0,0)=\lbrace (x,y) : x \not = 0, y = 0 \rbrace$. The Lagrangian descriptor defined in \eqref{M_function} for this system takes the form,
\begin{equation}
\begin{array}{ccl}
  M_p((x_0,y_0),t_0,\tau) & = & \displaystyle{\int^{t_0+\tau}_{t_0-\tau}} \left| x_0e^{F(t)} f(t) \right|^p +
\left| -y_0e^{-F(t)} f(t) \right|^p dt = |x_0|^p A(t) + |y_0|^p B(t)
\end{array}
\label{eq:lagrangian_non_aut_saddle_point}
\end{equation}
where
\begin{equation}
A(t) = \int^{t_0+\tau}_{t_0-\tau} \left| e^{F(t)} f(t) \right|^{p} \; dt \quad,\quad B(t) = \int^{t_0+\tau}_{t_0-\tau} \left| e^{-F(t)} f(t) \right|^{p} \; dt
\end{equation}

\noindent are functions that do not depend on $x_0$ and $y_0$. Consequently, \eqref{eq:lagrangian_non_aut_saddle_point} has the same functional form as \eqref{eq:lagrangian_aut_saddle_point} and we can apply the same argument as given in Theorem
\ref{thm:lin_aut_saddle_point} to obtain non-differentiability of $M_p$ along any line transverse to the stable manifold $W^{s}(0,0)$} for $p \leq 1$. Similarly we obtain non-differentiability of $M_p$ along any line transverse to the unstable manifold $W^{u}(0,0)$.

\subsection{The nonautonomous nonlinear saddle point}
\label{sec:nonlinear_nonauto_saddle}

We now consider a nonautonomous nonlinear system having a hyperbolic saddle trajectory at the origin. Let $f \in C^1(t_0-\tau, t_0+\tau)$, and $f(t) > 0$ for every $t \in [t_0-\tau,t_0+\tau]$. Our system has the form,
\begin{equation}
\left \{ \begin{array}{ccl}
   \dot{x} & = & f(t) x  + g_1(t, x, y) \\
   \dot{y} & = & -f(t) y + g_2(t, x, y), \\
\end{array}\right .
\label{eq:nonlin_non_aut_saddle_point}
\end{equation}

\noindent where $g_1, g_2: \mathbb{R} \times \mathbb{R}^2 \rightarrow \mathbb{R}$ with $g_1(t,0,0)=g_2(t,0,0)=0$, $\forall t$. We suppose that $g_1, g_2$ are real valued nonlinear functions and their order is quadratic or higher in $(x, y)$, with $g_1, g_2 \in C^1$ satisfying the conditions that make \eqref{eq:nonlin_non_aut_saddle_point} Hamiltonian. It is straightforward to verify that $(x, y) = (0, 0)$ is a hyperbolic trajectory. For hyperbolic trajectories in nonautonomous vector fields, their corresponding stable and unstable manifolds are, respectively, tangent to the stable and unstable subspaces of the linear approximation at the hyperbolic trajectory (\cite{CL, irwin, deblasi, KH, fen91}). We will show that the Lagrangian descriptor detects the manifolds in this example in the same manner as for the earlier examples in this paper. 

The strategy for demonstrating this fact is exactly the same as the one used to show that the structure of the Lagrangian descriptor for the linear autonomous saddle (discussed in Section \ref{sec:linear_auto_saddle}) coincides with the structure for the autonomous nonlinear saddle point (discussed in Section \ref{sec:nonlinear_auto_saddle}). A 
differentiable, invertible change of coordinates was made (given by \cite{moser56}), in such a way that effectively transformed the autonomous nonlinear saddle point into the form of the autonomous linear saddle point. Then it was obvious that the 
results for the autonomous linear saddle point, in the transformed coordinates, carried over the autonomous nonlinear saddle point directly. A final argument showing that the singularities of the Lagrangian descriptor for the nonlinear autonomous saddle point in the transformed coordinates were also present in the original coordinates completed 
the argument.

A similar strategy is carried out in the nonautonomous case, but it is not as straightforward as in the autonomous case. There has been a recent extension of Moser's theorem to the nonautonomous case (\cite{FW15}), but the requirements on the time dependence are too stringent for all the applications that we will consider, and therefore we will not utilize this result in our arguments for this example. Another approach is to utilize  a result like the Hartman-Grobman theorem (\cite{hart60a,hart60b,hart63,grob59,grob62}) for nonautonomous systems. Recently this result has been generalized to nonautonomous systems in \cite{bv06}. However, an issue with these ``linearization theorems'' is that  the linearization transformations are not differentiable on the hyperbolic trajectory. This situation has been discussed in detail in \cite{carlos} where it is argued that, for our purpose, this situation is essentially of a technical nature and does not prevent our use of such results to reach the desired conclusion.

\section{Non-Hamiltonian systems}
\label{sec:non_Ham_sys}

In this section we consider some issues related to the interpretation  of the output of Lagrangian descriptors when they are applied to non-Hamiltonian systems.

As an example, we consider the following non-Hamiltonian system:
\begin{equation}
\begin{cases}
\dot{x} = \lambda x \\
\dot{y} = -\mu y \\
\end{cases}\;,\quad \mu , \lambda>0 \,\, {\rm and}\,\, \mu \neq \lambda \label{eqds}
\end{equation}

\noindent for which the exact solution takes the expression,
\begin{equation}
\begin{cases}
   x(t) = x_0 e^{\lambda t}\\
   y(t) = y_0 e^{-\mu t}\\
\end{cases}
\end{equation}

\noindent In this case the origin is a hyperbolic fixed point with stable and unstable manifolds as in the example described in Section \ref{sec:linear_auto_saddle},
\begin{equation}
\begin{split}
  W^s(0,0) = \left\lbrace (x,y) \in \mathbb{R}^2 \;|\; x=0, y \neq 0 \right\rbrace, \\
  W^u(0,0) = \left\lbrace (x,y) \in \mathbb{R}^2 \;|\; y=0, x \neq 0 \right\rbrace.
\end{split}
\end{equation}

\noindent For the dynamical system (\ref{eqds}) at $t_0 = 0$, $M_p((x_0,y_0),t_{0},\tau) $ is explicitly computed as:
\begin{equation}
  \begin{array}{ccl}
    M_p((x_0,y_0),t_0,\tau) & = & \displaystyle{\int^{\tau}_{-\tau} |\lambda x_0 e^{\lambda t}|^p + |-\mu
y_0 e^{-\mu t}|^p dt } \\
      &   &                \\
      & = & \displaystyle{\lambda^{p-1}|x_0|^p \frac{(e^{\lambda p \tau} - e^{-\lambda p \tau})}{p} + \mu^{p-1}| y_0|^p
\frac{(e^{\mu p \tau}- e^{-\mu p \tau})}{ p}}. \\
  \end{array} \label{mex}
\end{equation}

\noindent From the form of the LD it is easy to see that non-differentiability of the directional derivatives on the stable and unstable manifolds follows from the same argument given in Theorem \ref{thm:lin_aut_saddle_point}. The singularities in expression (\ref{mex}) satisfy the Definition \ref{def}, which in turn is stated in the spirit of the phenomenology described  in \cite{prl,cnsns}. Therefore it is useful to point out here
that the term ``singularities''  does {\bf not} refer to any 
property of the contour lines of $M_p$, as incorrectly asserted in  (\cite{RH2015}).  Based on this misinterpretation \cite{RH16arxiv} has claimed  that (\ref{eqds}) is  a `'counter-example'' to the method of Lagrangian descriptors. This is an incorrect   statement as   demonstrated by (\ref{mex}). 

\begin{figure}
\centering
a)\includegraphics[scale = 0.4]{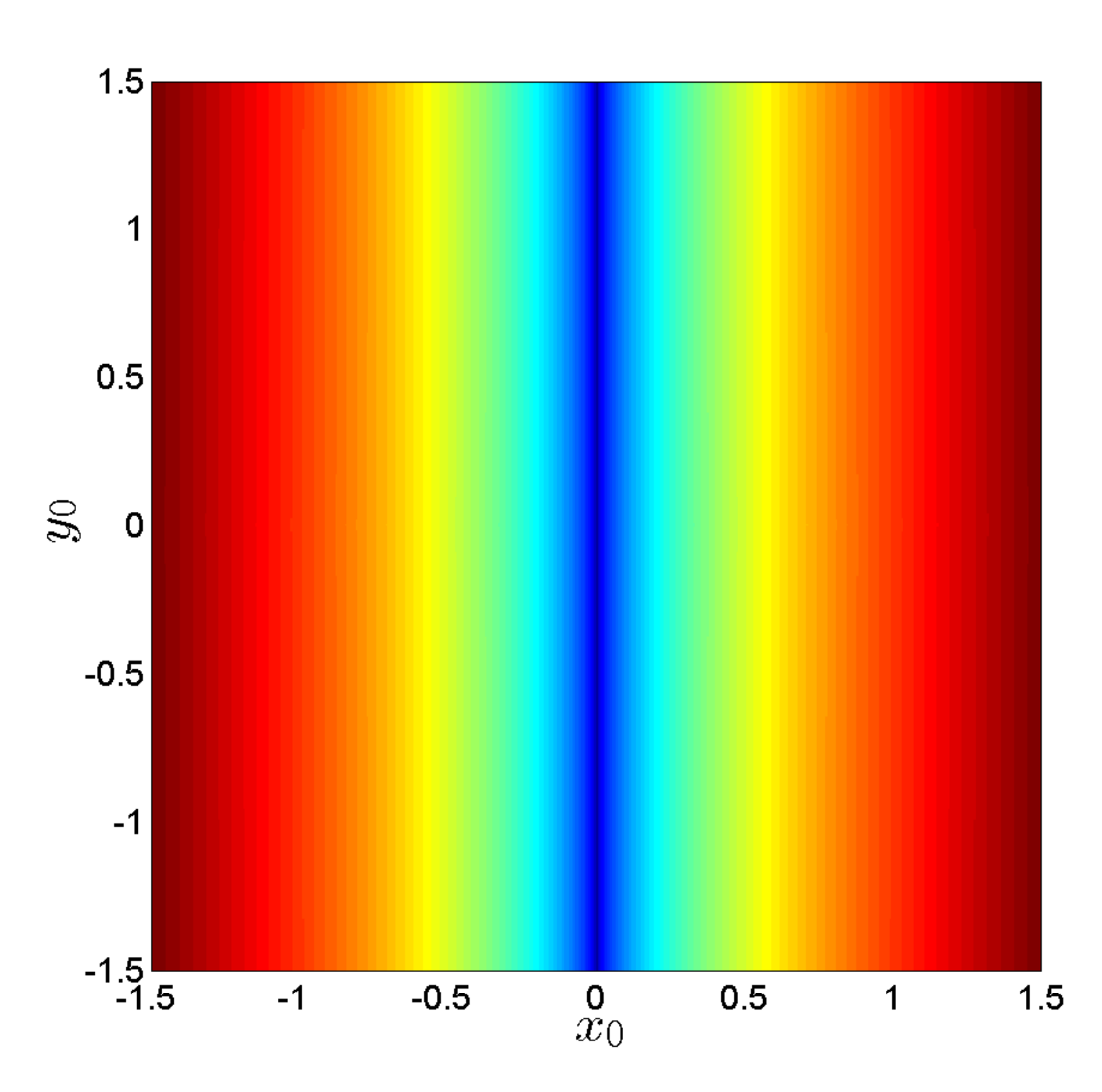}
b)\includegraphics[scale = 0.4]{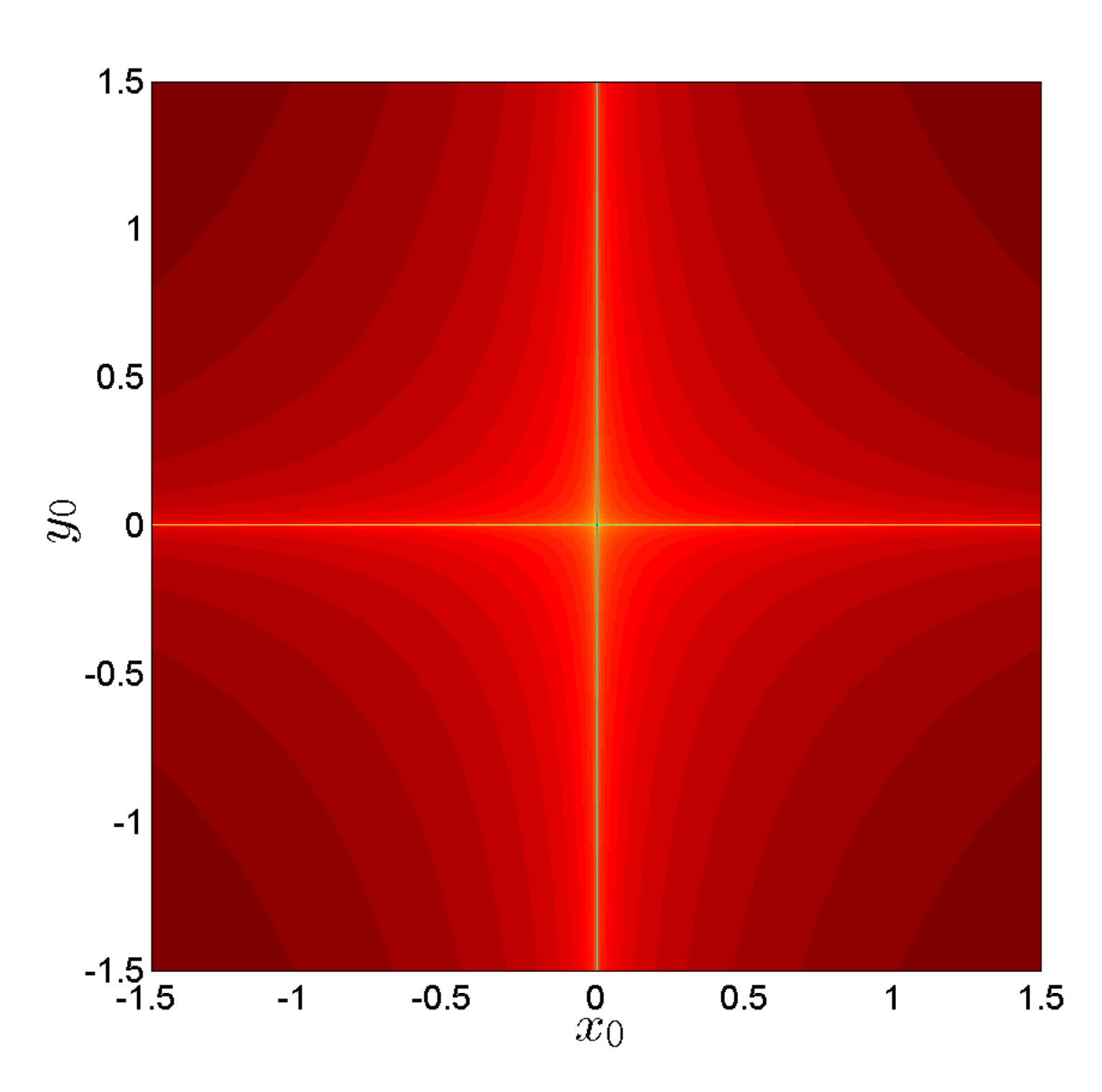}
\caption{Contour plots of Lagrangian descriptors computed for system (\ref{eqds}) with $\lambda = 2$ and $\mu = 1$: a) $M_p$ for $p = 0.5$ and $\tau=15$; b) $M_p$ for $p = 1/\tau$ and $\tau=15$.}
\label{fig:contours_M_vs_Mp}
\end{figure}

\begin{figure}
\centering
a)\includegraphics[scale = 0.4]{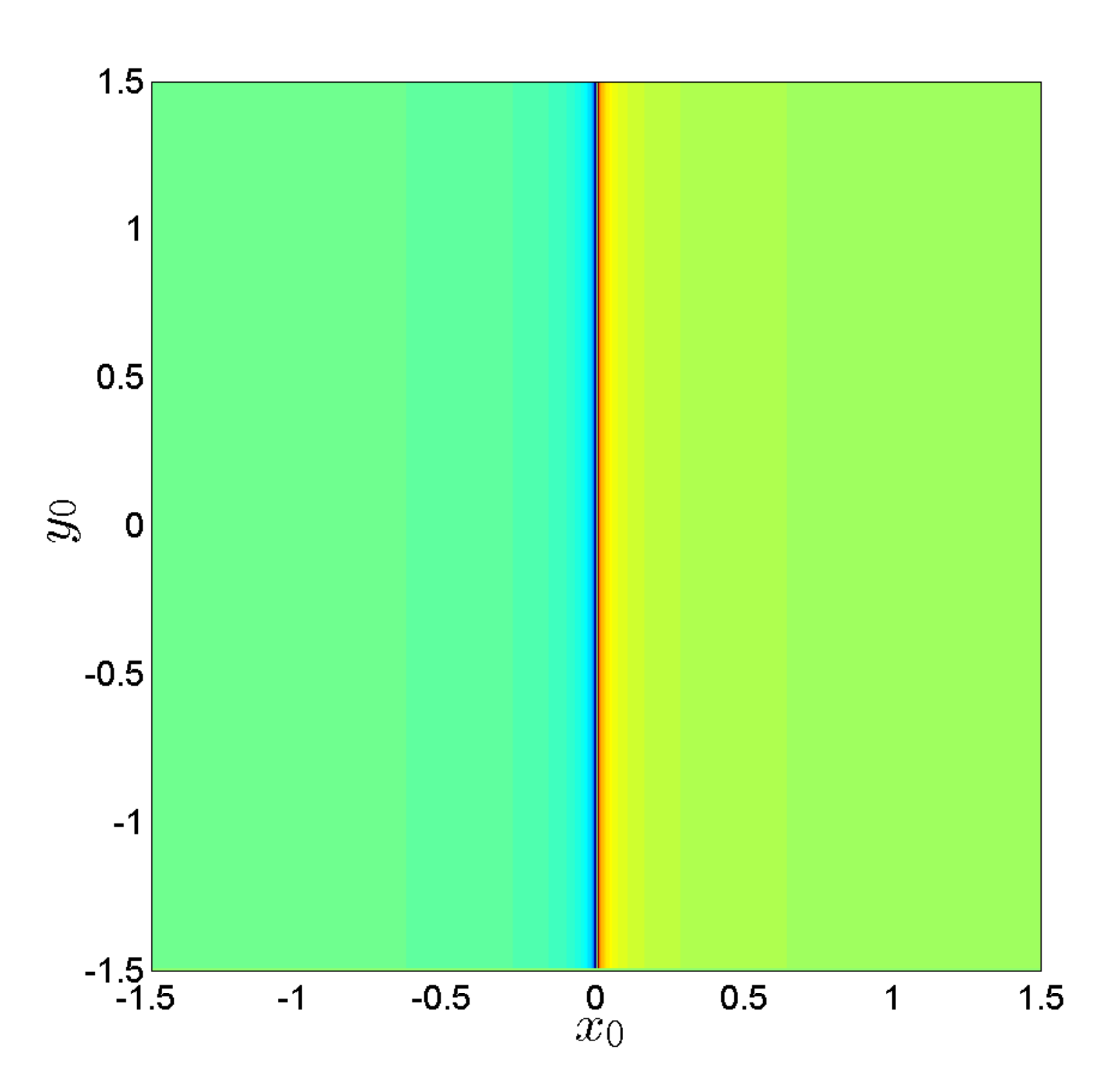}
b)\includegraphics[scale = 0.4]{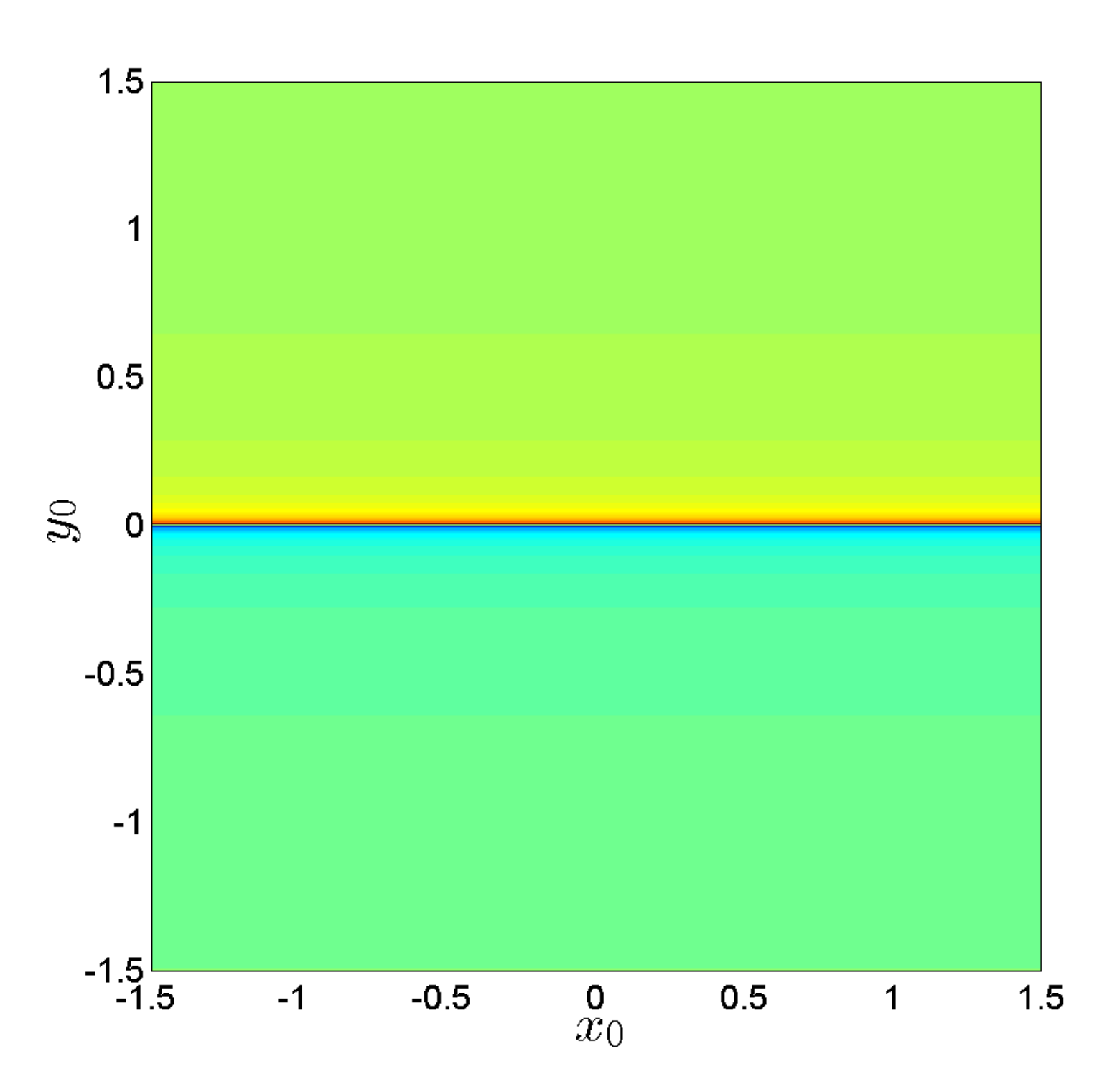}
\caption{Lagrangian descriptors applied to the system (\ref{eqds}) with $\lambda = 2$ and $\mu = 1$: a) Contours of $\frac{\partial M_{p}}{\partial x_{0}}$ using $p=0.5$ and $\tau=15$; b) Contours of $\frac{\partial M_{p}}{\partial y_{0}}$ using $p=0.5$ and $\tau=15$.}
\label{fig:partial_derivative_M_p}
\end{figure}

This implies that  care should be taken when trying to visualize the singular features of (\ref{mex}) as in Definition \ref{def} from contour plots of $M_{p}$ itself. Observe that given a fixed $p$ and a sufficiently large $\tau$, the terms in expression (\ref{mex}) can take values which can differ by orders of
magnitude. Figure \ref{fig:contours_M_vs_Mp} a) illustrates this point by showing the 
contours of $M_p$ with $p = 0.5$ and $\tau=15$, for the values $\lambda = 2$ and $\mu = 1$. In this case the first term in (\ref{mex}) is much greater than the second, and thus the effect of the  latter goes unnoticed in the figure. Still, the presence 
of the singularities can also be highlighted by plotting the partial derivatives of $M_p$, as illustrated in Figure  \ref{fig:partial_derivative_M_p}. 
Using the quotient of the two terms in expression (\ref{mex}), where each term in the quotient evaluates  $M_p$ on each manifold, we obtain an idea of  their different contribution.
\begin{equation}
\displaystyle{\frac{M_p((x_0,0),0,\tau)}{M_p((0,y_0),0,\tau)} = \left (\frac{\lambda}{\mu} \right )^{p-1} \frac{|x_0|^p}{|y_0|^p} \frac{e^{2p\tau\lambda}-1}{e^{2p\tau\mu}-1}
e^{p\tau(\mu - \lambda)}}.\nonumber
\end{equation}

\noindent In particular, for large $\tau$ we have
\begin{equation}
\displaystyle{\frac{M_p((x_0,0),0,\tau)}{M_p((0,y_0),0,\tau)} = \left (\frac{\lambda}{\mu} \right )^{p-1} \frac{|x_0|^p}{|y_0|^p} e^{p\tau(\lambda - \mu)}}. \label{weight}
\end{equation}

\noindent Observe that a good visualization of both manifolds in the same plot can be obtained from the contours if the exponential in (\ref{weight}) is kept of order ${\cal O}(1)$ (see Fig. \ref{fig:contours_M_vs_Mp} b)). This is achieved for $p\tau(\lambda - \mu) \approx 1$, that is,
\begin{equation}\label{eq:p_formula}
p = \frac{1}{\tau(\lambda - \mu)}.
\end{equation}

We consider next another non-Hamiltonian system,
\begin{equation}
\begin{cases}
\dot{x} = - x, \\
\dot{y} = - y. \\
\end{cases}
\label{eqglobalatt}
\end{equation}

\begin{figure}[H]
\centering
a)\includegraphics[scale = 0.4]{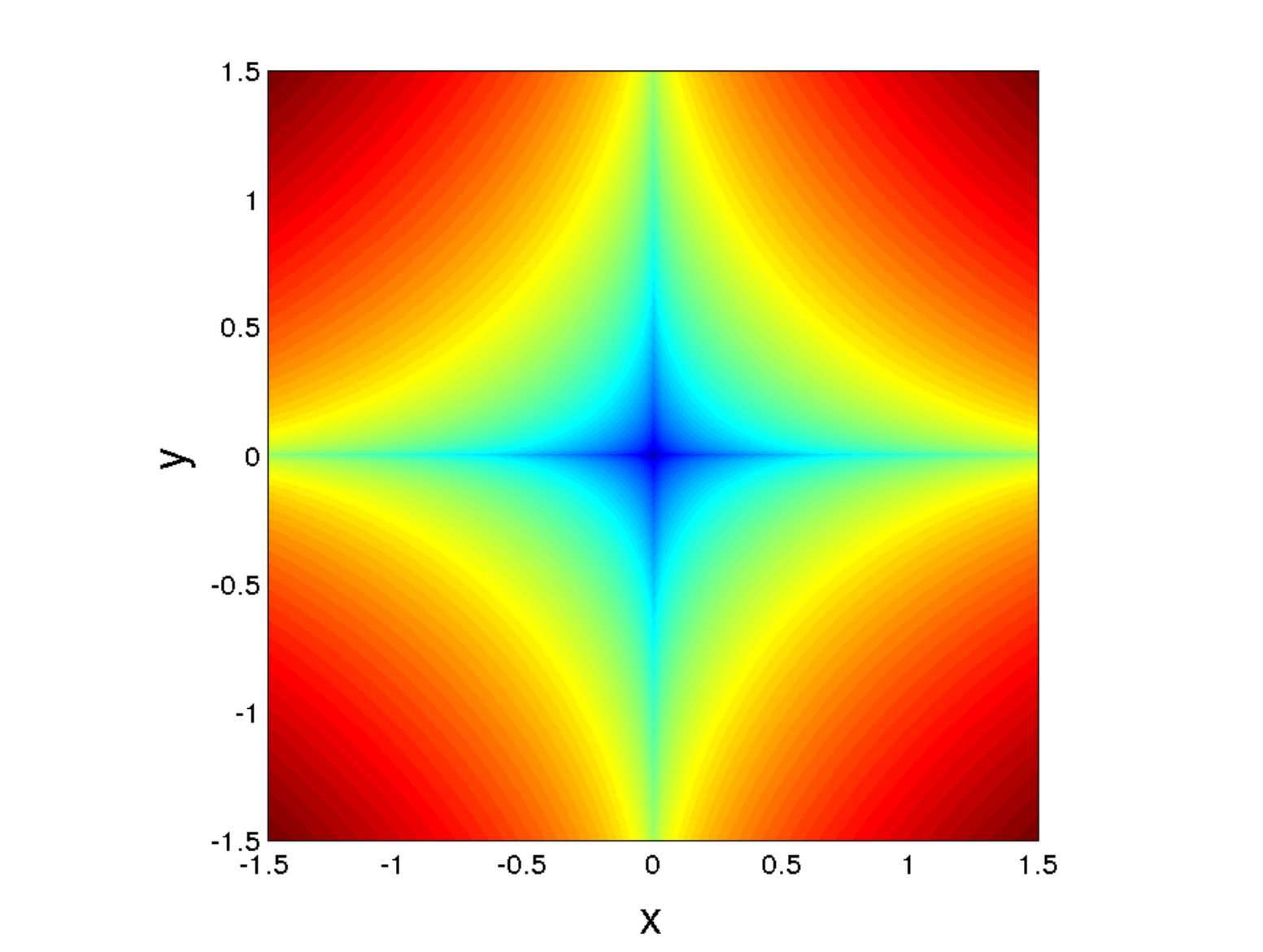}
b)\includegraphics[scale = 0.4]{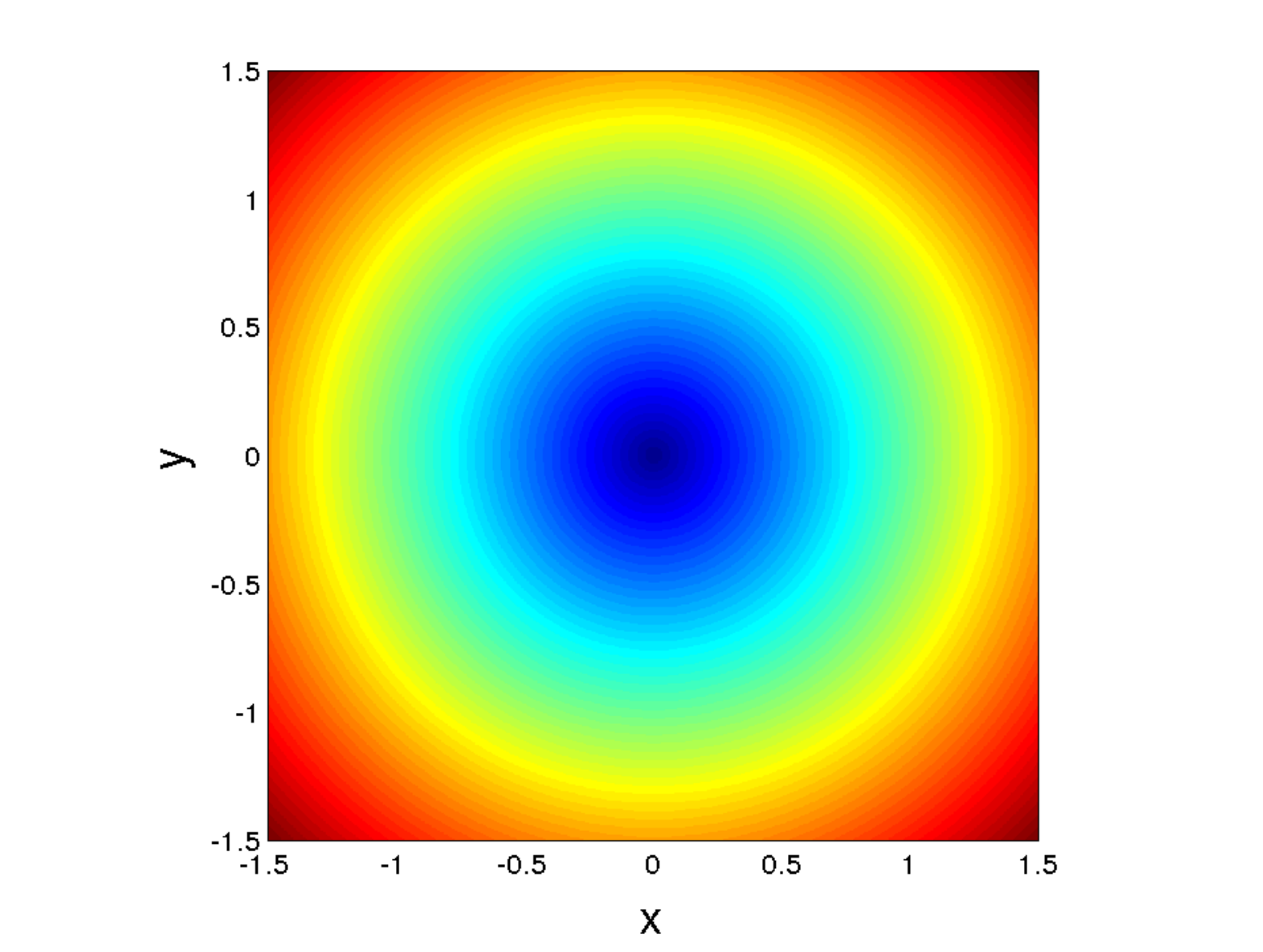}
\caption{The left-hand panel shows contours of $M_p$ for $p=0.5$ and the right-hand panel shows contours of $M$ 
function. Both are obtained for system \eqref{eqglobalatt} for $\tau=15$.}
\label{fig:globalatt}
\end{figure}

\noindent This system has a global attractor at the origin. Furthermore, we can analytically compute Lagrangian descriptors for this example. In the case of the $M_p$ function, we have
\begin{equation}
  \displaystyle{M_p((x_0,y_0),0,\tau) = \int^{\tau}_{-\tau} |x_0 e^{-t}|^p + |y_0 e^{-t}|^p dt = 
(|x_0|^p+|y_0|^p)\frac{e^{-p\tau} + e^{p\tau}}{p}}
\end{equation}

\noindent where the singularities are located on the axes $x=0$ and $y=0$, as is observed in Figure \ref{fig:globalatt}a). Certainly these lines correspond to stable manifolds of the fixed point, however any line passing through the origin  is a  stable manifold  of it, because in fact the whole plane is a stable manifold. In this case  one could ask about to what extent it is useful   highlighting just two lines of the whole plane. 

In addition, this observed feature at representing the $M_{p}$ function for the global attractor is not longer reproduced when computing original $M$ function as in \cite{chaos}. We then obtain
\begin{equation}
  \displaystyle{M((x_0,y_0),0,\tau) = \int^{\tau}_{-\tau} \sqrt{x^2_0 e^{-2t} + y^2_0 e^{-2t}} dt = 2\sqrt{x^2_0 + y^2_0} \sinh(\tau).}
\end{equation}

\noindent In this case $M$ does not highlight any singular feature, except for the fixed point at the origin, and therefore the entire stable manifold is not associated with singular feature, as can be seen in Figure \ref{fig:globalatt}b).

\section{Lagrangian Descriptors and the Ergodic Partition theory}
\label{sec:LD_elliptic}

This section  illustrates, by using a very simple example, the link between Lagrangian descriptors and previous rigorous results on the ergodic partition theory (\cite{mezic3, susuki}).

We consider the dynamical system:
\begin{equation}
\begin{cases}
\dot{x} = y, \\
\dot{y} = -x \\
\end{cases}
\label{sys2} 
\end{equation}

\noindent where the origin is an elliptic fixed point. This system can be expressed in action-angle variables $(\rho, \theta)$ as follows,
\begin{equation}
\begin{cases}
\dot{\rho} = \frac{\partial H}{\partial \theta} = 0 \\[.1cm]
\dot{\theta} = -\frac{\partial H}{\partial \rho} = -1
\end{cases}
,\quad \rho \in [0,\infty) \;,\; \theta \in [0,2\pi)
\label{sys3}
\end{equation}

\noindent where $x= \rho \cos \theta$ and $y= \rho \sin \theta$. The Hamiltonian expressed in action-angle variables is $H(\rho, \theta)=\rho$. From these expressions it is clear that $(\rho=\rho_0, \theta(t)=-t+\theta_0)$ are solutions to the system (\ref{sys3}) which correspond to  invariant 1-tori. In summary, the phase plane of system (\ref{sys2}) is foliated by invariant sets consisting of concentric circles.  
The autonomous system (\ref{sys2}) does not have hyperbolic fixed points nor their invariant manifolds, thus the ideas based on ``singular features'' of $M_p$ explained in the preceding sections are not directly applicable here. However it is interesting to realize that the described results of LDs are complementary to other previous results in the literature, which we explain next, that rigorously justify the ability of LDs to highlight Lagrangian coherent structures characterized by invariant $n$-tori.

The definition of $M_p$ given in (\ref{M_function}) is the sum of $n$ quantities representing the integration along trajectories of the functions $|\dot{x}_i|^p$, where $i \in \{1,\ldots,n\}$. These quantities, up to a factor $1/(2\tau)$, are the time averages along trajectories of the given functions. Analogously, the original definition of LDs given by \cite{prl,cnsns} based on the Euclidean arc length is obtained by integrating the modulus of the velocity along trajectories, i.e.,
\begin{equation}
M(\mathbf{x}_0,t_0,\tau) = \int_{t_0-\tau}^{t_0+\tau} \left\| {\textbf{v}}(\mathbf{x}(t;\mathbf{x}_0),t) \right\| \; dt.
\end{equation}

\noindent
Again, $M$ is the average along trajectories of the function $\left\| {\textbf{v}}\right\|$ up to a factor $1/(2\tau)$. The role of averages of functions along trajectories for obtaining invariant sets is dicussed in works by \cite{mezic3, susuki}. There, the Birkhoff ergodic theorem is used. This theorem states that in the limit ${\tau \to \infty}$, averages of functions along trajectories of dynamical systems which preserve smooth measures  and are defined on compact sets do exist. Level sets of these limit functions are invariant sets. However, we note that the Birkhoff ergodic theorem has {\em not} been generalized to the case of aperiodically time-dependent vector fields.

We examine these ideas in the case of $M_p$ for the example (\ref{sys2}). First we remark that Hamiltonian dynamical systems preserve smooth measures as they are volume preserving. We show next that the limit of the time average can be analytically calculated for this system using $p = 1$. Considering the solutions of (\ref{sys2}), $M_{p=1}$ is
\begin{equation}
\label{mp1}
\begin{split}
M_{p=1}(\mathbf{x}_{0},0,\tau) &= \rho_0 \int_{-\tau}^{\tau} | \sin \left( -t + \theta_0 \right) | + | \cos \left( -t + \theta_0 \right) | \; dt \\
 &= \rho_0  \int_{-\tau+ \theta_0}^{ \theta_0} |\sin(s)| + |\cos(s)| \; ds  + \rho_0 \int_{ \theta_0}^{\tau+ \theta_0} |\sin(s)| + |\cos(s)| \; ds .
\end{split}
\end{equation}

\noindent Now we can write $\tau = N\pi + q$ with $q \in [0,\pi)$ and calculate one of the integrals of (\ref{mp1}) as,
\[
\begin{split}
\int_{\theta_0}^{\theta_0+\tau} |\sin s| \; ds & = \int_{\theta_0}^{\theta_0+N\pi+q} |\sin s| \; ds = \int_{\theta_0}^{\theta_0+N\pi} |\sin s| \; ds + \int_{\theta_0+N\pi}^{\theta_0+N\pi+q} |\sin s| \; ds \\ & = N \int_{\theta_0}^{\theta_0+\pi} |\sin s| \; ds + \int_{\theta_0}^{\theta_0+q} |\sin s| \; ds = 2N + \int_{\theta_0}^{\theta_0+q} |\sin s| \; ds .
\end{split}
\]

\noindent Analogously for the other terms we obtain,
\[
\int_{\theta_0}^{\theta_0+\tau} |\cos s| \; ds = 2N + \int_{\theta_0}^{\theta_0+q} |\cos s| \; ds ,
\]
\[
\int_{-\tau+\theta_0}^{\theta_0} |\sin s| \; ds = 2N + \int_{\theta_0-q}^{\theta_0} |\sin s| \; ds ,
\]
\[
\int_{-\tau+\theta_0}^{\theta_0} |\cos s| \; ds = 2N + \int_{\theta_0-q}^{\theta_0} |\cos s| \; ds .
\]
Consequently the time average of $M_{p=1}$ is the limit\,
\begin{eqnarray}
\lim_{\tau \rightarrow \infty} \frac{1}{2\tau} M_1(\mathbf{x}_{0},0,\tau) &=&\lim_{N \rightarrow \infty}  \frac{8N\rho_0}{2(N\pi+q)} +  \frac{\rho_0}{2 \left( N\pi+q \right) }\int_{\theta_0-q}^{\theta_0+q} |\cos s|+|\sin s| \; ds\nonumber\\
&=&\frac{4\rho_0}{\pi}= \frac{4}{\pi} \sqrt{x_0^2+y_0^2}.
\end{eqnarray}

\noindent In order to prove that the time average of $M_p$ also converges for $0<p<1$ we can use that
$$\begin{array}{ccccl} 0 & \leq & \displaystyle{\frac{1}{2\tau} M_p(\mathbf{x}_{0},0,\tau)} & = & \displaystyle{\frac{1}{2\tau} \rho_0^p \int_{-\tau}^{\tau}| \sin \left( -t + \theta_0 \right) |^{p} + | \cos \left( -t + \theta_0 \right) |^{p} \; dt} \\
 & \leq & \displaystyle{\frac{1}{2\tau} \rho_0^p \int_{-\tau}^{\tau} (1+1) \; dt} & = & 2\rho_0^p ,\quad \quad \forall p \in (0,1].
\end{array}$$

\noindent Since $\frac{1}{2\tau} M_p(\mathbf{x}_{0},0,\tau)$ is an increasing function of $\tau$ and it is bounded by a constant value, then $\frac{1}{2\tau} M_p(\mathbf{x}_{0},0,\tau)$ also converges when $\tau \rightarrow \infty$.

\begin{figure}
\centering
a){\includegraphics[scale = 0.42]{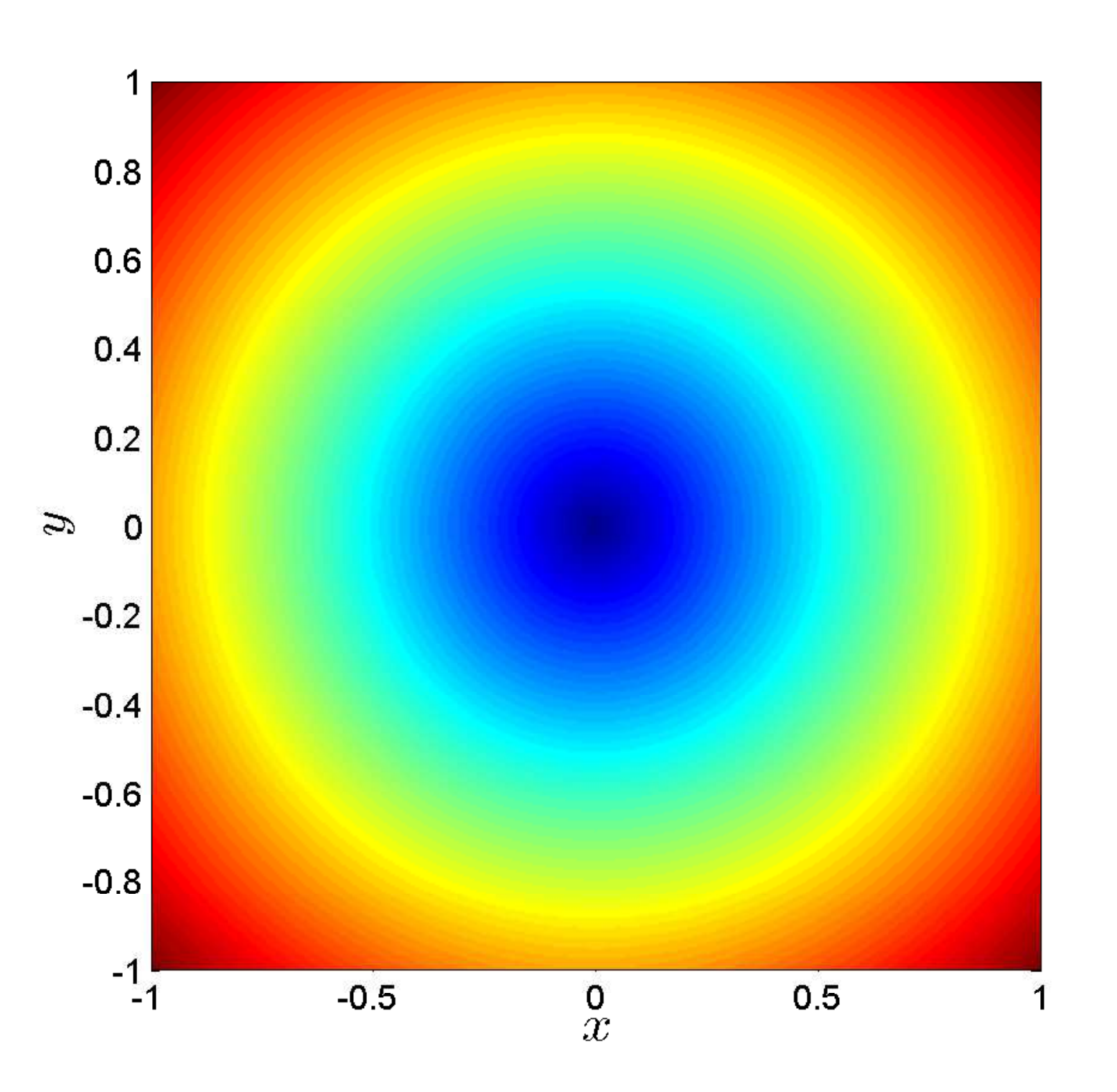}}
b){\includegraphics[scale = 0.42]{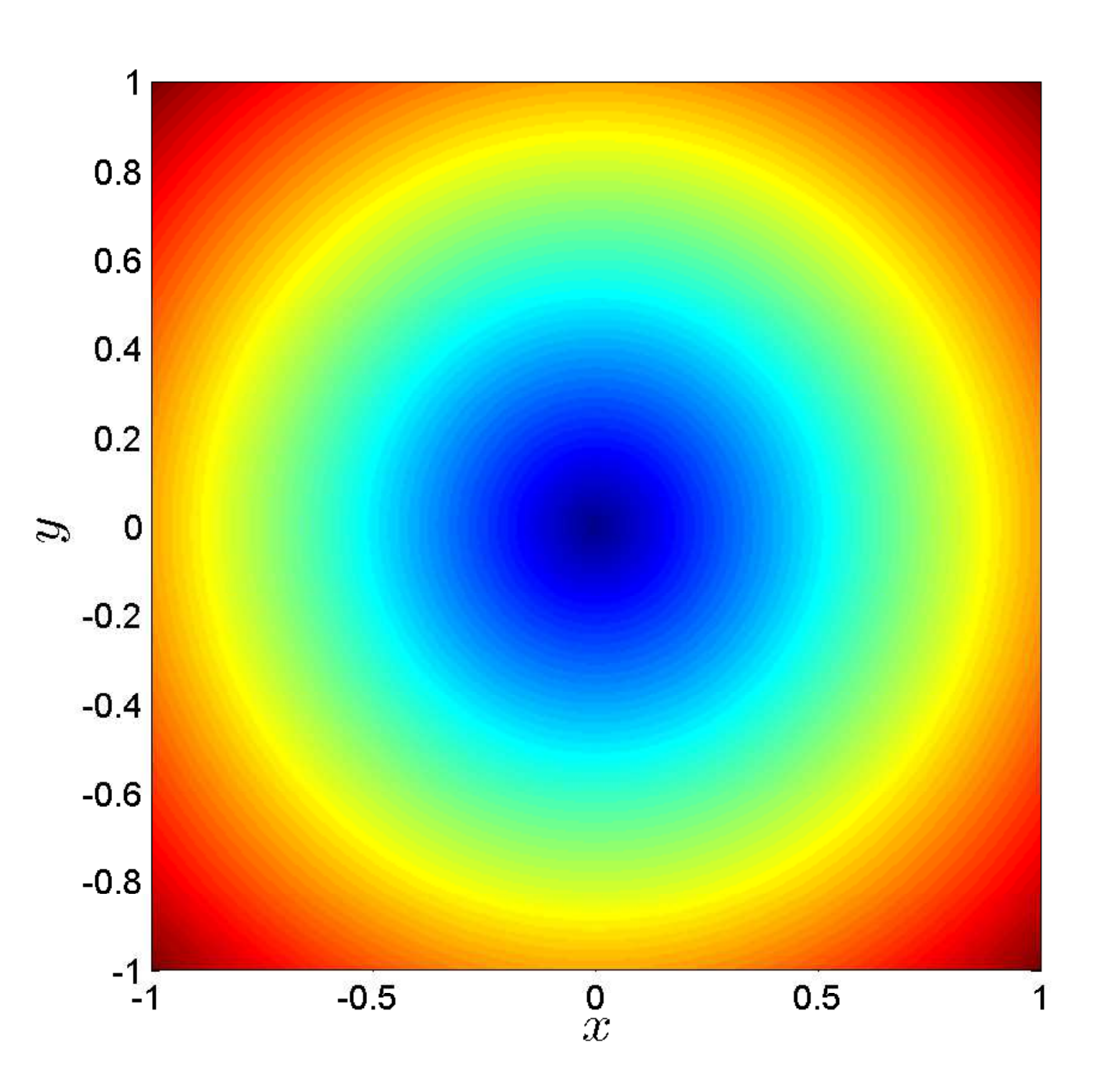}}
\caption{Phase space of (\ref{sys2}) calculated with: a) $M$ for $\tau = 10$; b) $M_1$ for $\tau = 10$.}
\label{fig:M_tori}
\end{figure}

Results for $M$ are also easily obtained, given that arc length is $M = 2\tau\rho$. In order to obtain the limit of the time average here, it is not necessary to go up to very large $\tau$, since the time average is a constant function in $\tau$ (i.e. $\frac{1}{2\tau}M = \rho$), and thus the convergence is obtained for any finite $\tau$. We note that contour lines of either $M_p$ and $M$ are the same as their averages. The important point for those level sets to be invariant sets is that they have to be taken once the convergence of the average is reached. In the case of $M_p$ a sufficient large $\tau$ is required and in the case of $M$ any $\tau$ is valid. In both cases the invariant sets are just the concentric circles given by the 1-tori (see Figure \ref{fig:M_tori}). Further averages of functions to complete an ergodic partition as described by \cite{mezic3, susuki}, are not required for this particular example, as the circles are already minimal invariant sets.
 
As a final remark we note that despite the similarity between Figures \ref{fig:M_tori}a) and  \ref{fig:globalatt}b),  the interpretation of LDs is completely different  in each case. The link  between contour lines of $M$  and invariant sets  requires a measure preserving  dynamical system defined on a compact set, and this is clearly not the case for system (\ref{eqglobalatt}) represented in  Fig. \ref{fig:globalatt}b).  
 
\section{Lagrangian Descriptors and 3D flows}
\label{sec:3D}

Visualizing flow structures in three dimensions is of much interest, but achieving success requires overcoming numerous difficulties. Probably the fundamental issue is that it is difficult to ``organize'' the data from an ensemble of trajectories in such a way as to ``reveal'' geometrical structures.  In this context, LDs were applied to study the transport in the three-dimensional unsteady Hill's spherical vortex in \cite{cnsns}.  The method of Lagrangian descriptors was successful in this study in that it revealed  both invariant manifolds of hyperbolic trajectories and invariant sets related to $n$-tori solutions.

A brief review on the background and issues associated with the dynamical systems approach to transport in three dimensions was given in \cite{wiggins2010}. A collection of representative references of the dynamical systems approach to Lagrangian transport in three dimensions, that should not be interpreted as an exhaustive review of this topic, are:  \cite{mackay,mw1,fountain,svl2001,xiumei,mjm2005,green07,bw09,psc2010,mezic10, lsmr2012,slm2012,slrm2014,McIlhany,srlm2016}.

We illustrate the ability of LDs to visualize three dimensional flows by analyzing the well known Arnold-Beltrami-Childress (ABC) flow (\cite{arnold65,arnold83,dombre}). This flow models prototypes of fast dynamos in magnetohydrodynamics (cf. \cite{arnold83,galloway}) and it is also a steady solution of Euler's equations for inviscid fluid flows (see \cite{child}). 

The equations for fluid particle trajectories of the ABC flow are given by\textcolor{magenta}{,}
\begin{equation}
\begin{cases}
\dot{x} = A \sin z + C \cos y \\
\dot{y} = B \sin x + A \cos z \\
\dot{z} = C \sin y + B \cos x
\end{cases} ,\;\; x,y,z \in [0,2\pi] \;.
\label{ds_abc}
\end{equation} 

 \noindent 
where $A$, $B$, and $C$ are parameters to be chosen shortly. The ABC flow is one of the first flows for which the existence of chaotic particle paths was demonstrated (cf. \cite{arnold65}). Subsequently, there have been numerous studies of the flow structure of the ABC flow from the dynamical systems point of view, e.g., \cite{dombre,haller_physd,sulman}. Here we show that both the very intricate manifold structures and the coherent structures  of the ABC flow can be visualized with Lagrangian descriptors. In order to illustrate this we choose two sets of parameter values,
\[
A = B = C = 1 \;,
\]
\[
A = 1\;,\;B = \sqrt{\dfrac{2}{3}}\;,\;C = \dfrac{\sqrt{3}}{3} \;.
\]

\begin{figure}[htbp!]
  \centering
  \subfigure[$A = B = C = 1$]{\includegraphics[scale=0.53]{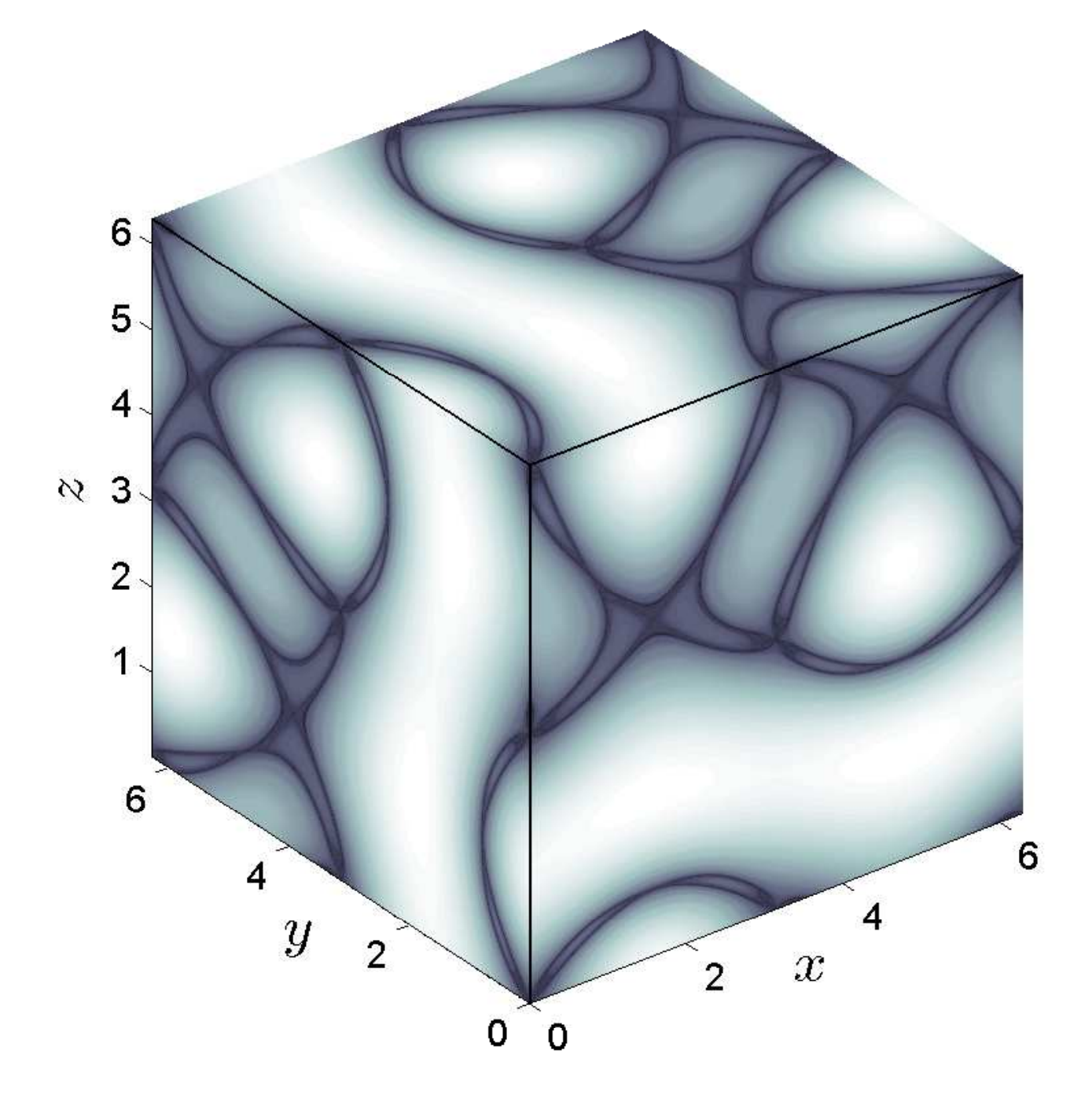}}
  \subfigure[$A = 1\;,\;B = \sqrt{\dfrac{2}{3}}\;,\;C = \dfrac{\sqrt{3}}{3}$]{\includegraphics[scale=0.53]{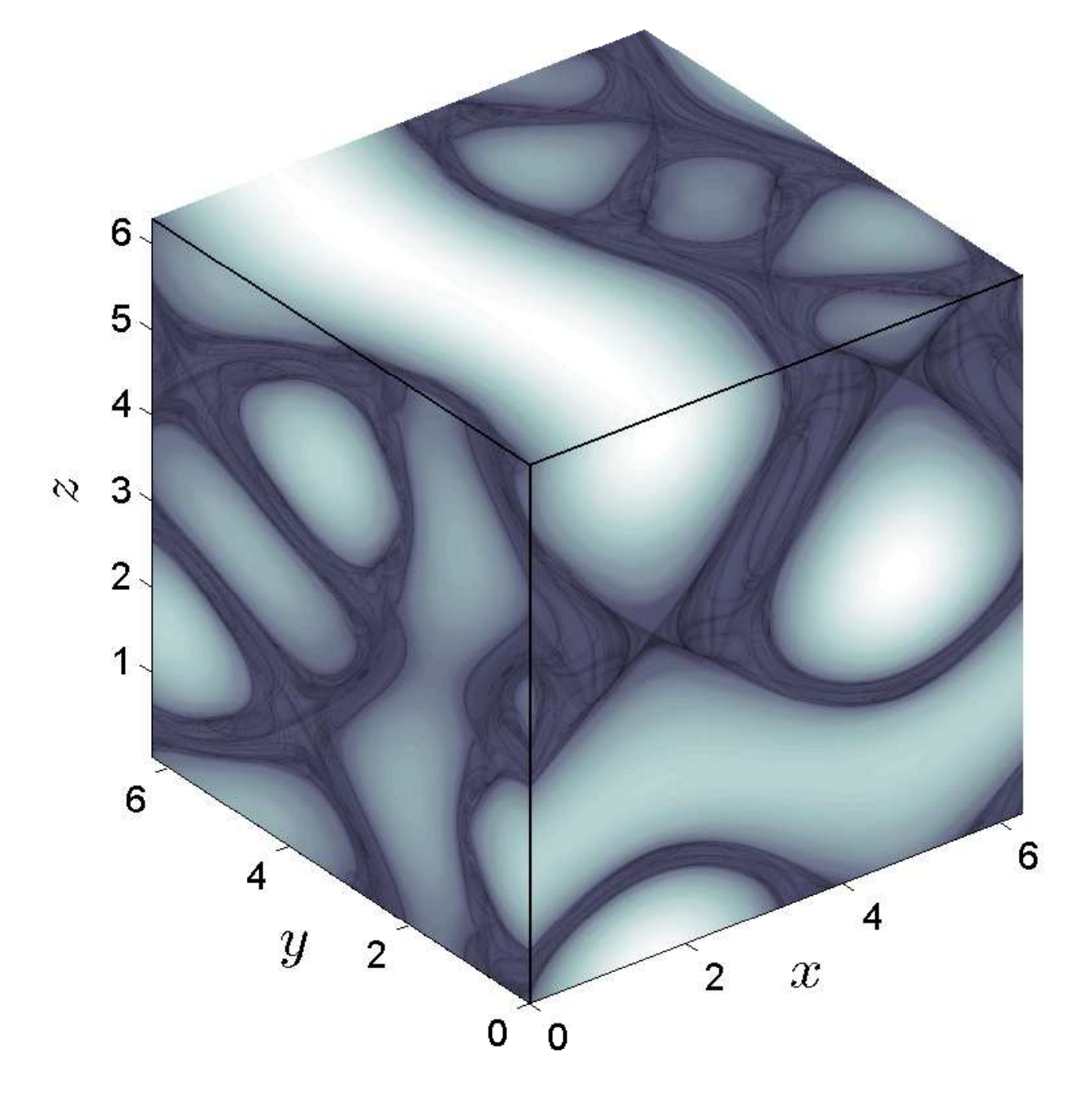}} \\
	\subfigure[$A = B = C = 1$]{\includegraphics[scale=0.53]{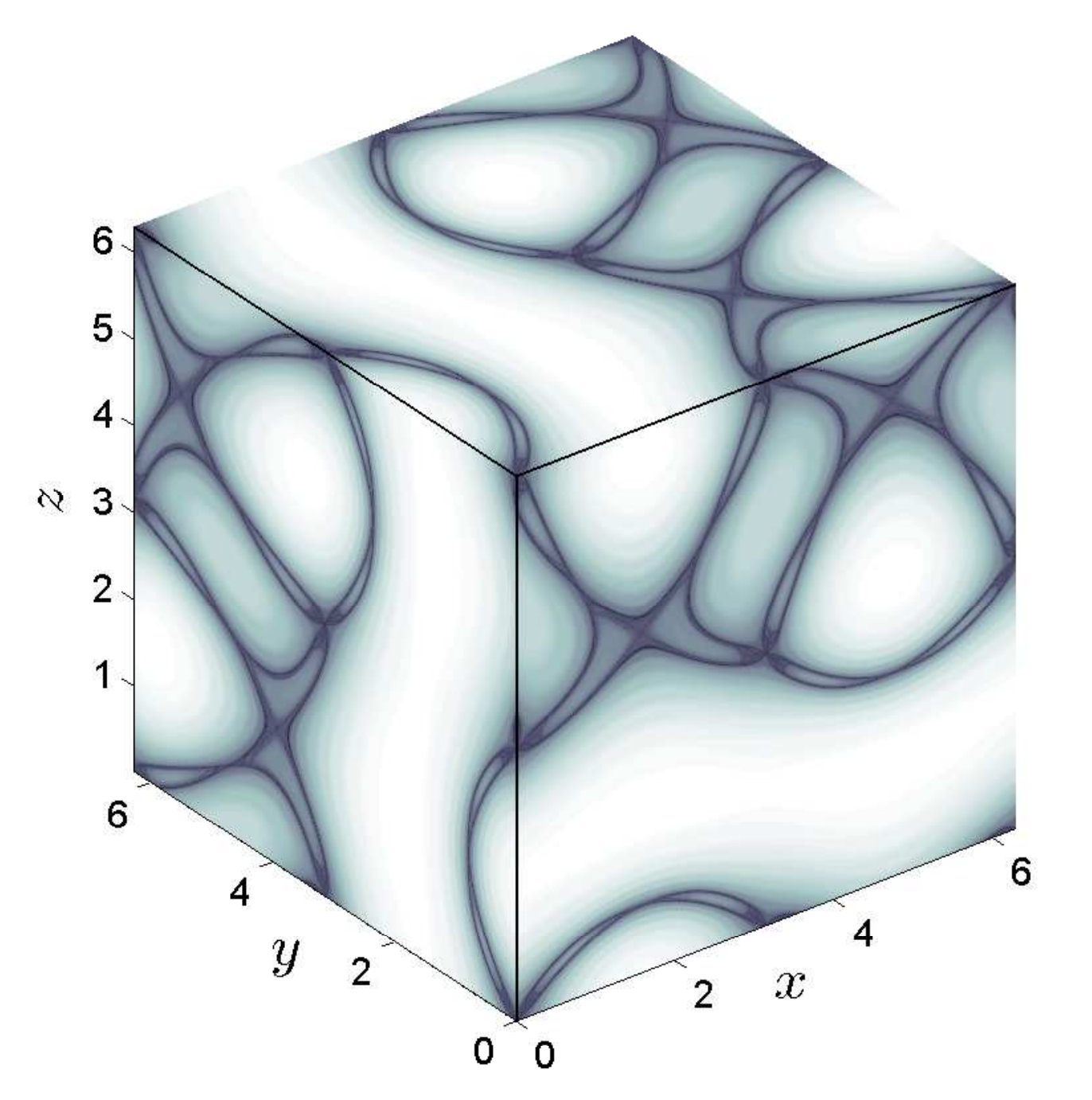}}
  \subfigure[$A = 1\;,\;B = \sqrt{\dfrac{2}{3}}\;,\;C = \dfrac{\sqrt{3}}{3}$]{\includegraphics[scale=0.53]{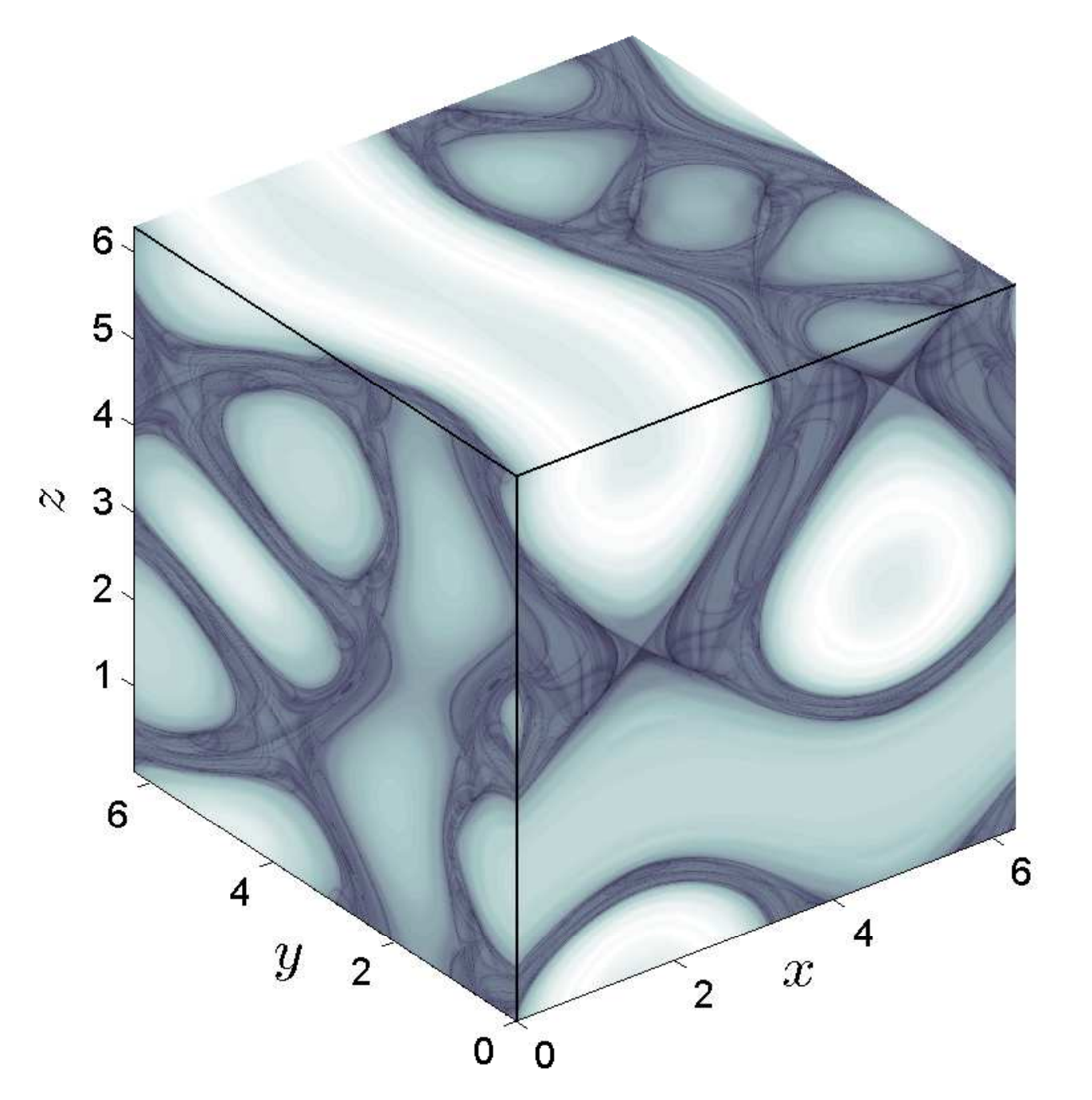}}
  \caption{A) and B) Contours of $M$ for the ABC flow with $\tau = 30$; C) and D) Contours of $M_1$ for the ABC flow with $\tau = 30$.}
  \label{ABC_3D}
\end{figure}

\begin{figure}[htbp!]
  \centering
  \subfigure[$A = B = C = 1$]{\includegraphics[scale=0.48]{fig8a.pdf}}
  \subfigure[$A = 1\;,\;B = \sqrt{\dfrac{2}{3}}\;,\;C = \dfrac{\sqrt{3}}{3}$]{\includegraphics[scale=0.53]{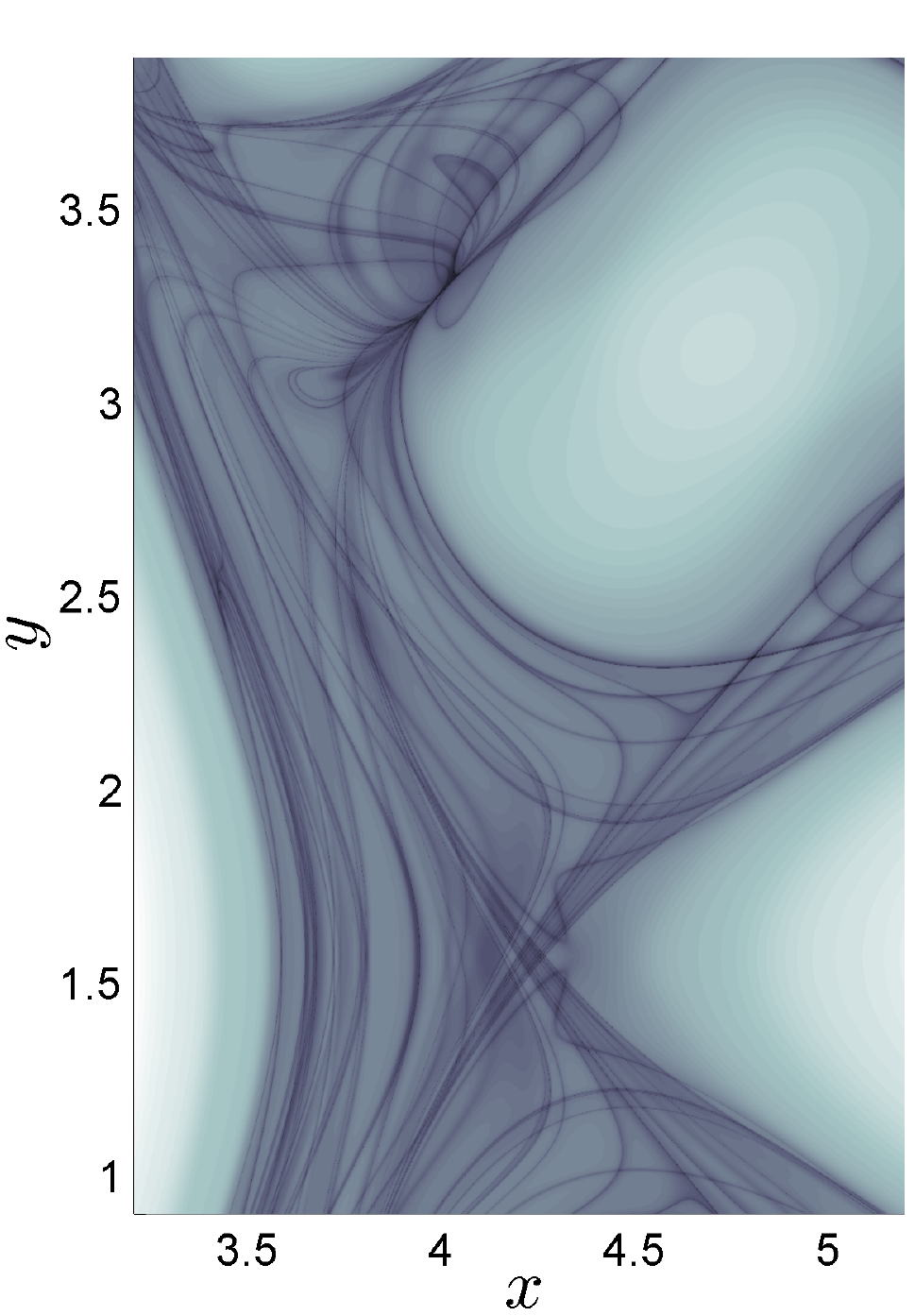}} \\
	\subfigure[$A = B = C = 1$]{\includegraphics[scale=0.48]{fig8c.pdf}}
  \subfigure[$A = 1\;,\;B = \sqrt{\dfrac{2}{3}}\;,\;C = \dfrac{\sqrt{3}}{3}$]{\includegraphics[scale=0.53]{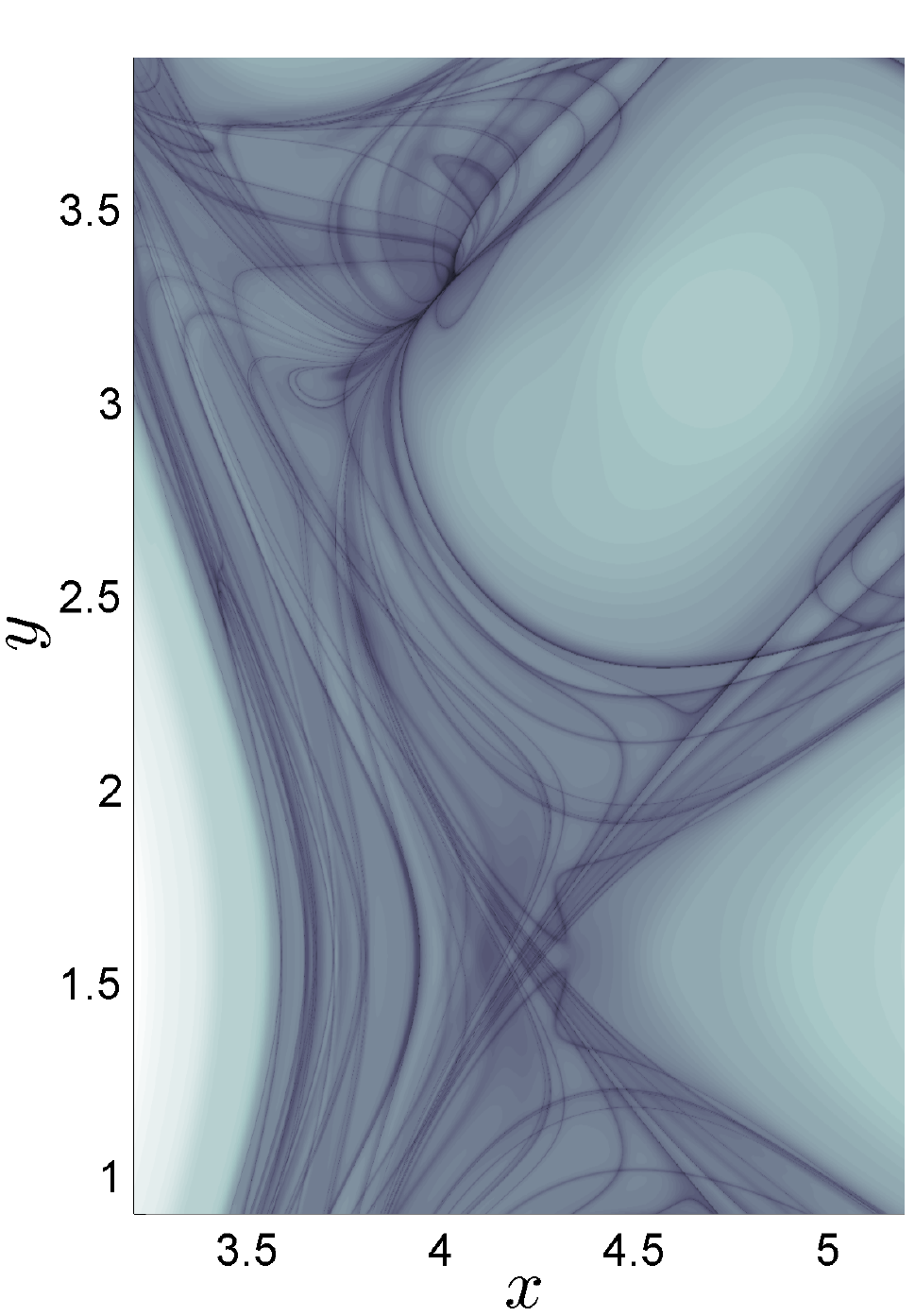}}
  \caption{Contours of LDs with $\tau = 30$: (a) zoom of $M$ on the plane $\{y=0\}$; (b) zoom of $M$ on the plane $\{z=0\}$; (c) zoom of $M_1$ on the plane $\{y=0\}$; d) zoom of $M_1$ on the plane $\{z=0\}$.}
  \label{ABC_plane}
\end{figure}

\noindent Figure \ref{ABC_3D} shows the results obtained with $M$ and $M_1$ for these parameter values. The figure shows both coherent structures and the complicated tangle of repeatedly intersecting stable and unstable manifolds of hyperbolic trajectories along heteroclinic orbits. This chaotic tangle provides the ``geometric  template'' for the chaotic mixing mechanism of particles, which is visible from Figure \ref{ABC_plane}.

In addition, we demonstrate next with this example the link of LDs to invariant sets. The results discussed in the previous section are applicable here as the ABC flow {is incompressible and consequently preserves an invariant measure. In order to analyze the detection of KAM tori by means of LDs for the ABC flow we will focus on the system with parameters $A = 1\;,\;B = \sqrt{2}/\sqrt{3}\;,\;C = \sqrt{3}/3$. In particular, we consider a line of initial conditions,
\begin{equation}
\mathbf{x}_0 = (0,3.2,z_0) \quad,\; z_0 \in [3.6,5.9] 
\label{ic_line}
\end{equation}

\noindent which crosses an elliptic region of phase space as shown in Figure \ref{init_conds_ellip_an} with a blue. For these initial conditions we study the convergence of the time averages of $M$ and $M_p$ so the invariant sets present in the elliptic regions of the phase portrait for the ABC flow can be recovered from the contours of LDs contained in these regions. Figures \ref{t_avg_conv_M} and \ref{t_avg_conv_Mp} show the evolution of the time averages of $M/(2\tau)$ and $M_1/(2\tau)$ respectively. Observe also that their convergence at two particular initial conditions, one inside the elliptic region (marked with a magenta cross) and one inside the chaotic tangle regime (green cross), have been highlighted (see figs. \ref{init_conds_ellip_an}, \ref{t_avg_conv_M} and \ref{t_avg_conv_Mp}). These figures emphasize how at initial conditions inside elliptic regions, the time averages $M/(2\tau )$ and $M_{1}/(2\tau )$ reach convergence for sufficiently large $\tau$, meanwhile for initial conditions located in hyperbolic regions, $M/(2\tau )$ and $M_{1}/(2\tau )$ do not seem to converge for that time period (see figs. \ref{t_avg_conv_M}, \ref{t_avg_conv_Mp} and \ref{t_avg_conv_chaos}). Thus contour lines of $M$ and $M_p$ at the $\tau$ values where convergence is met are in 1-1 correspondence to invariant set. Thus contrary to what is stated by \cite{faraz}, LDs distinguish coherent structures which are invariant and this is backed by specific mathematical results. Similar results to the ones discussed here are found in \cite{budi}. 
We remark that the decomposition achieved by the LD is not minimal in the domain    
$x,y,z \in [0,2\pi]$. As an example, this  is observed from the red level sets displayed in the top and bottom right corners in Fig.  \ref{init_conds_ellip_an}b), which would correspond to  invariant sets disconnected from the one inside the elliptic region. Fig. \ref{traj_inv_set} shows the integration of a trajectory starting from the initial condition $(0,3.2,4.1)$, marked with a magenta cross on one of the contour lines obtained from the time-average convergence (see fig. \ref{init_conds_ellip_an}), confirming the invariant and minimal character of the level-set restricted to the elliptic region. 

\begin{figure}[htbp!]
\centering
A){\includegraphics[scale = 0.42]{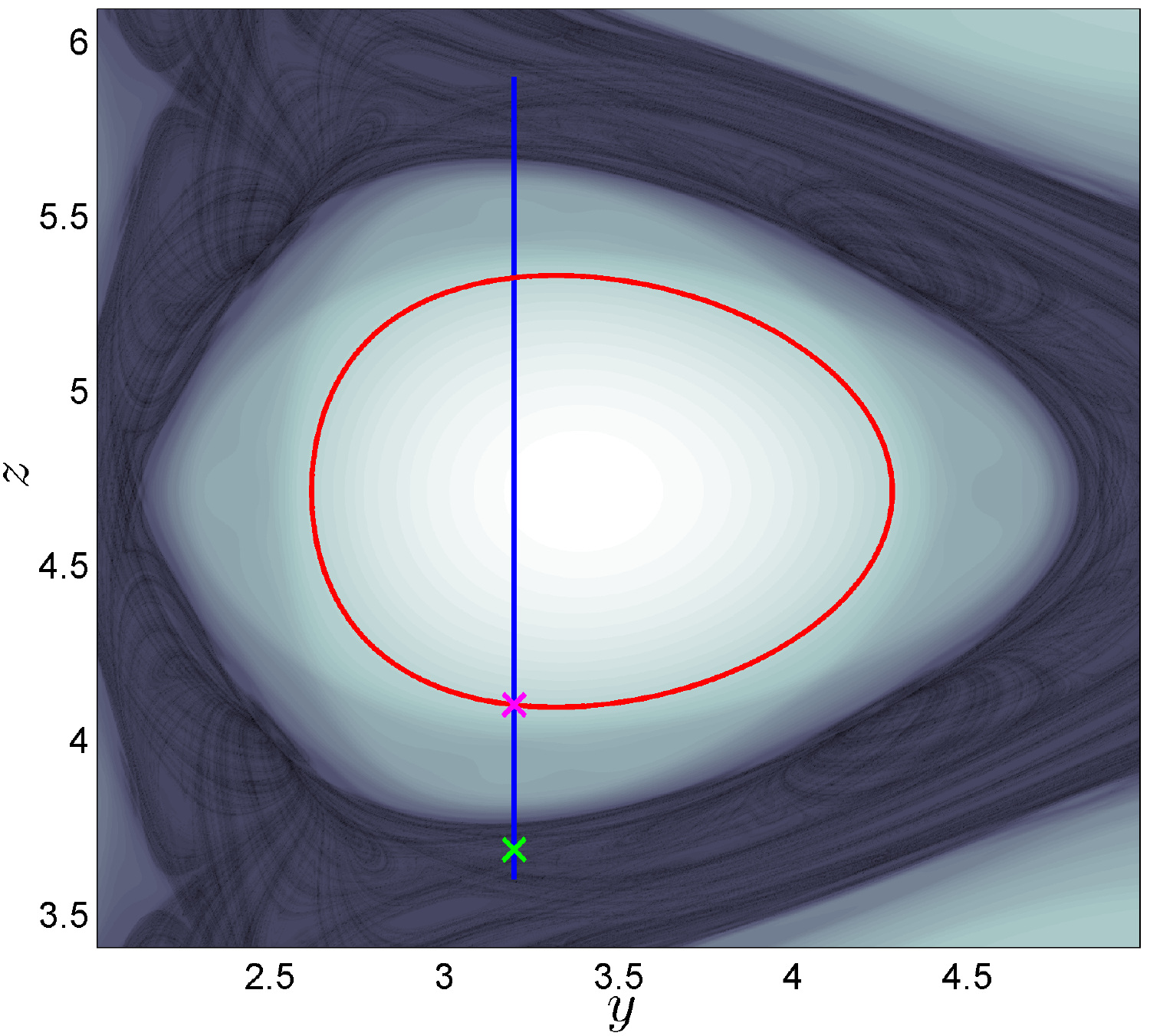}}
B){\includegraphics[scale = 0.42]{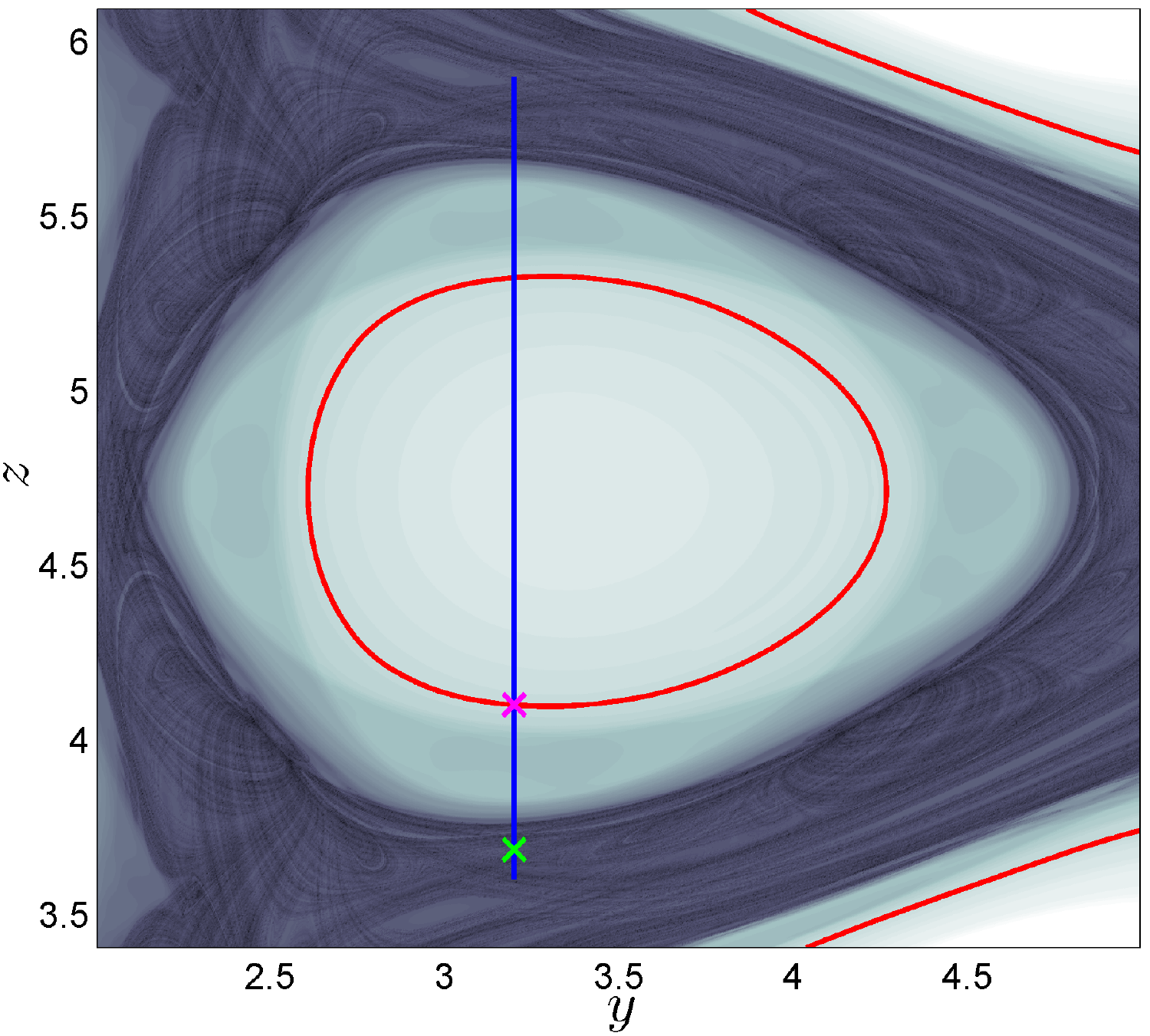}}
\caption{A) $M$ function on an elliptic region located on the plane $\{x = 0\}$ for $\tau = 75$ (value for which convergence inside the elliptic region is ensured); B) $M_1$ function on an elliptic region located on the plane $\{x = 0\}$ for $\tau = 100$ (value for which convergence inside the elliptic region is ensured). In both panels the blue line represents the line of initial conditions considered for the time average analysis, the magenta cross is an initial condition inside the elliptic region and the green cross an initial condition located in a chaotic region. Also, we display in red color the contours corresponding to the magenta initial condition after the time average has converged.}
\label{init_conds_ellip_an}
\end{figure}

\begin{figure}[htbp!]
\centering
a){\includegraphics[scale = 0.42]{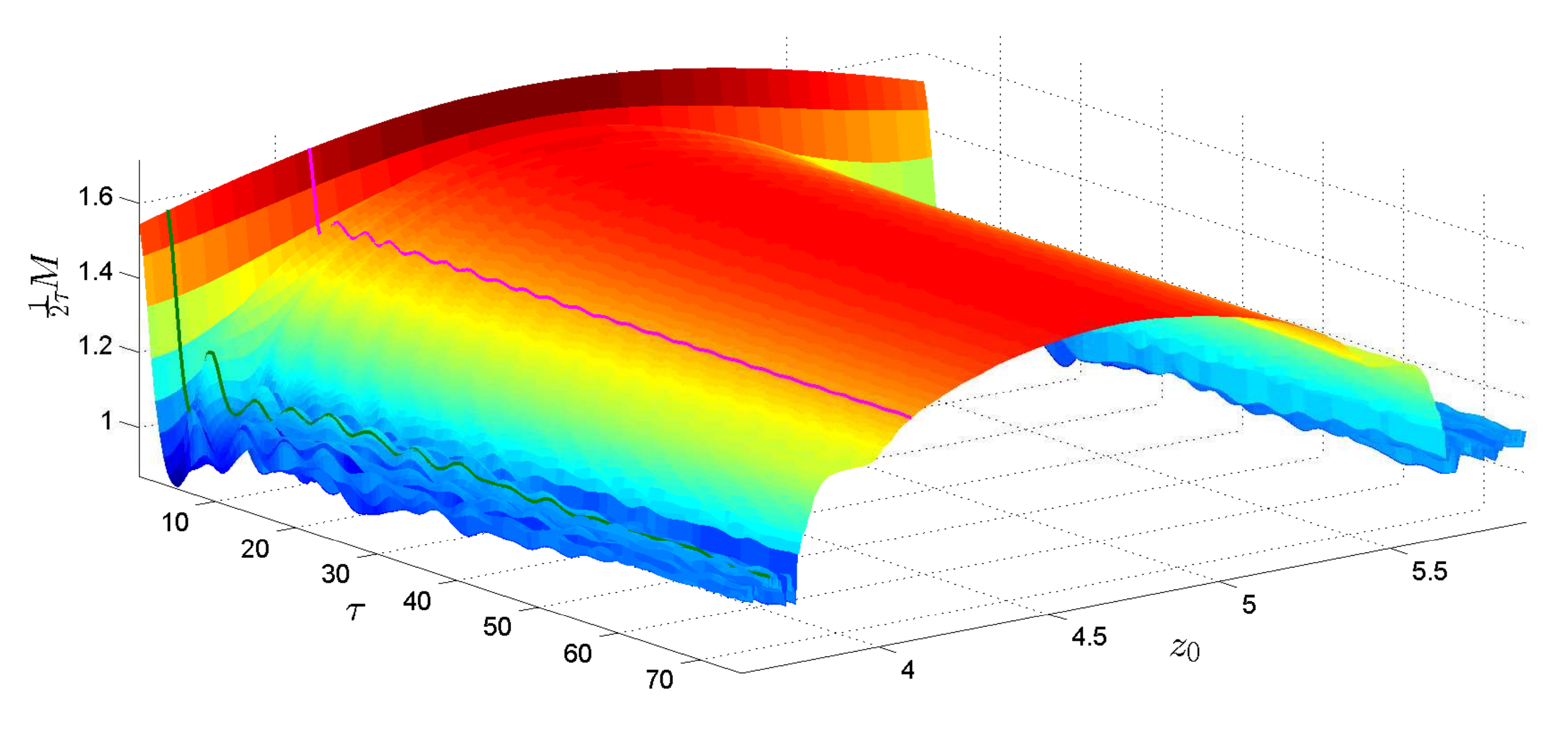}}
b){\includegraphics[scale = 0.42]{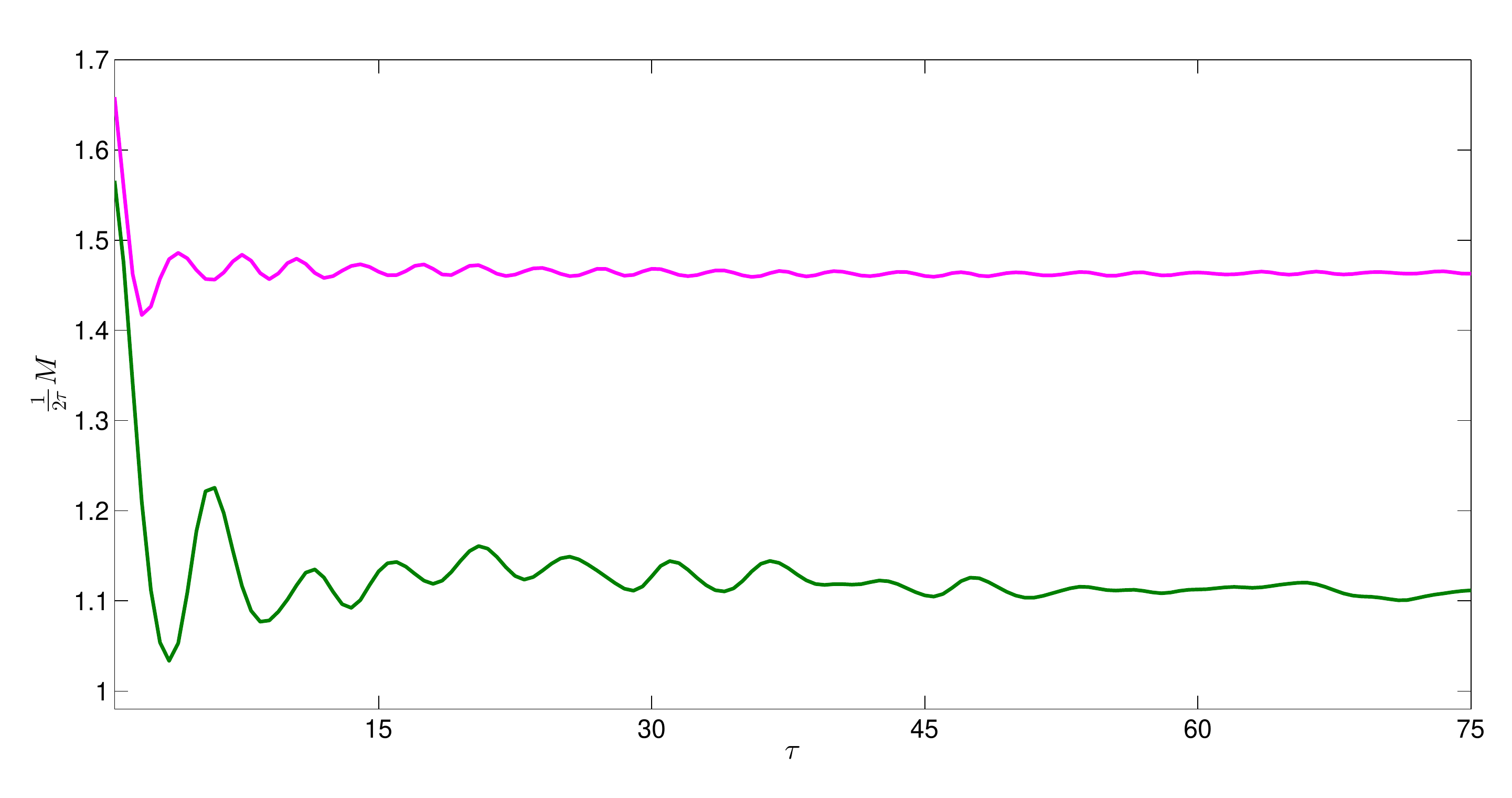}}
\caption{a) Time average evolution of $M$ in the range $\tau \in (0,75]$ for the line of initial conditions (\ref{ic_line}); b) Time average evolution of the pink and green initial conditions depicted in fig. \ref{init_conds_ellip_an}.}
\label{t_avg_conv_M}
\end{figure}

\begin{figure}[htbp!]
\centering
a){\includegraphics[scale = 0.42]{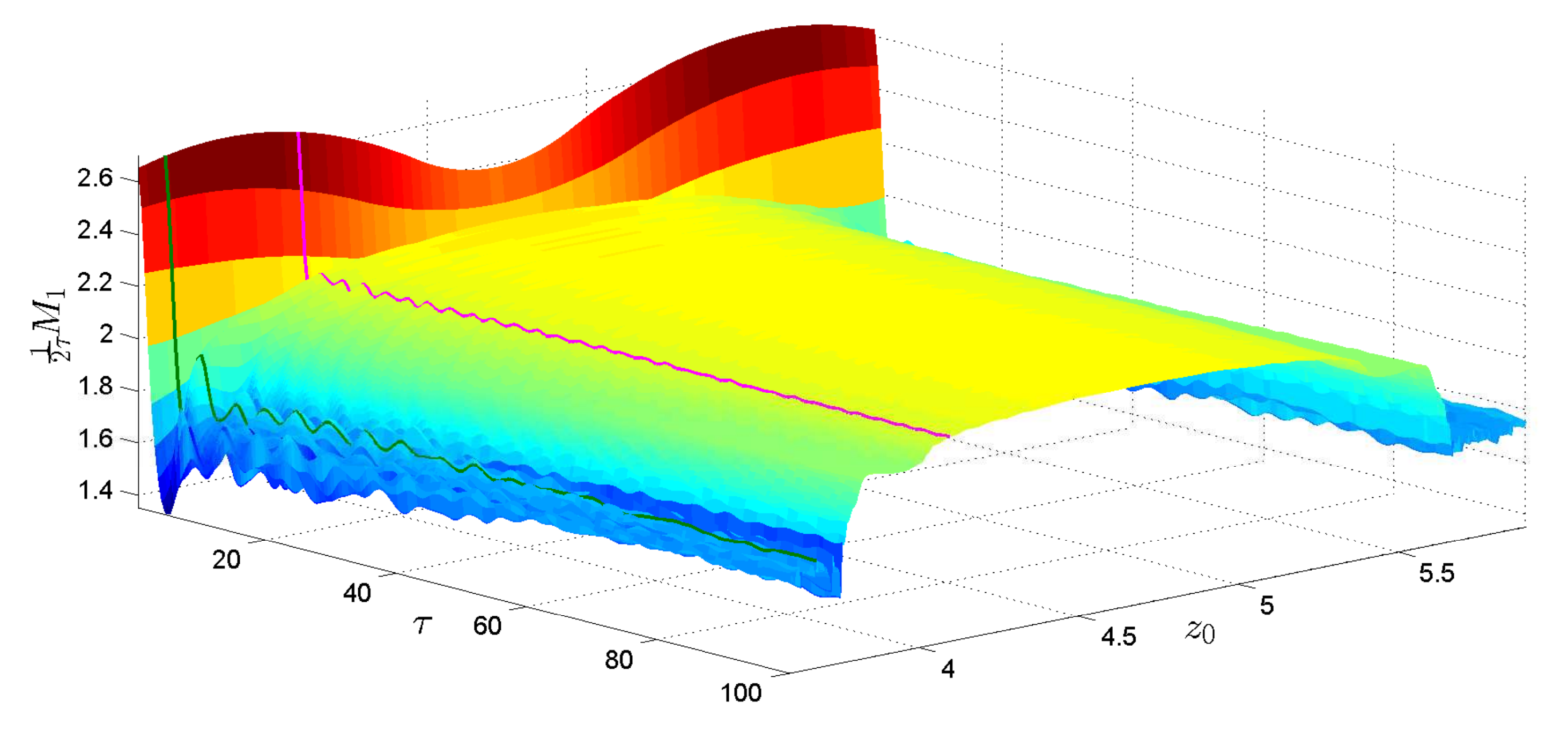}}
b){\includegraphics[scale = 0.42]{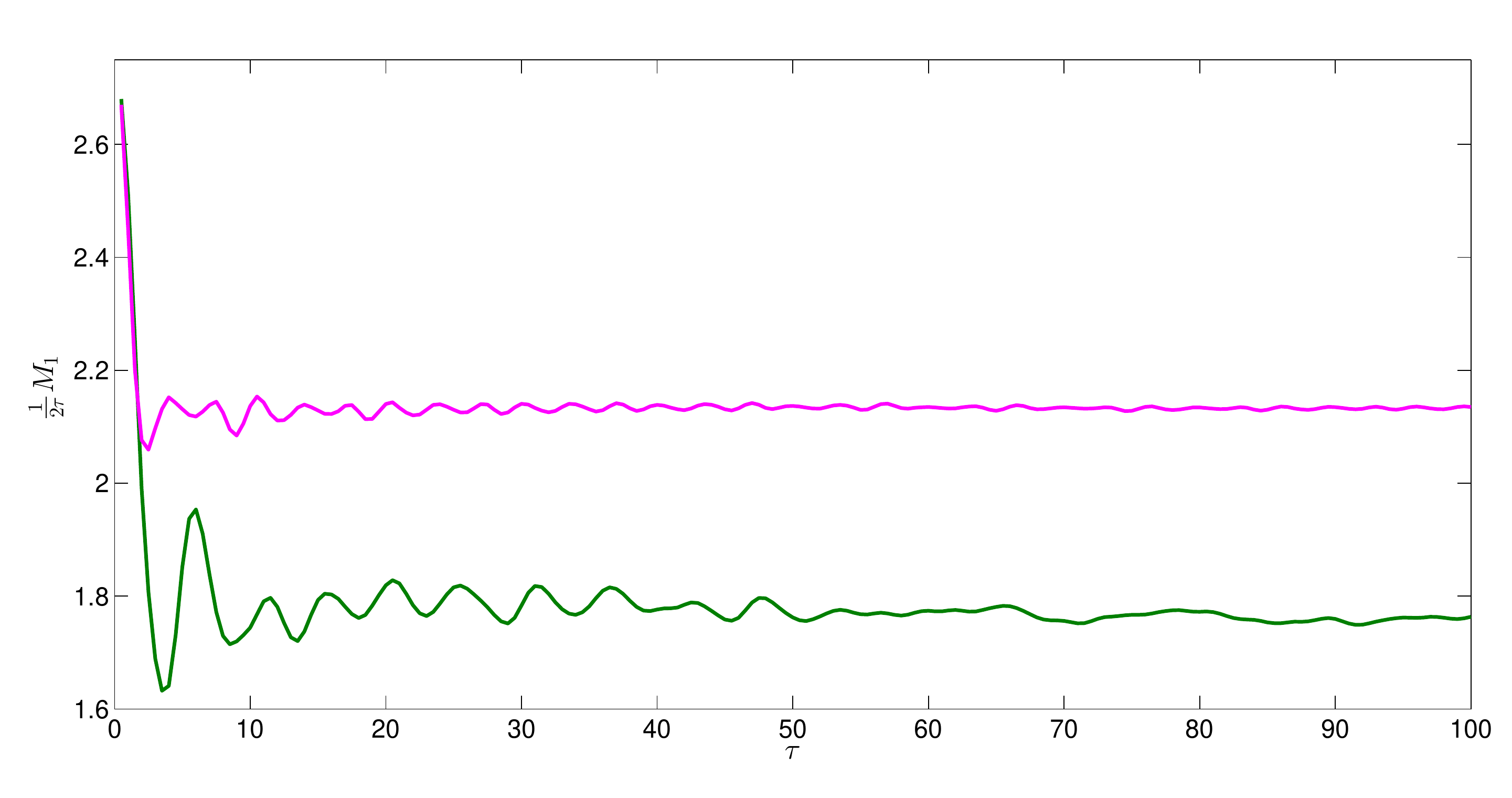}}
\caption{a) Time average evolution of $M_1$ in the range $\tau \in (0,100]$ for the line of initial conditions (\ref{ic_line}); b) Time average evolution of the pink and green initial conditions depicted in fig. \ref{init_conds_ellip_an}.}
\label{t_avg_conv_Mp}
\end{figure}

\begin{figure}[htbp!]
\centering
a){\includegraphics[scale = 0.42]{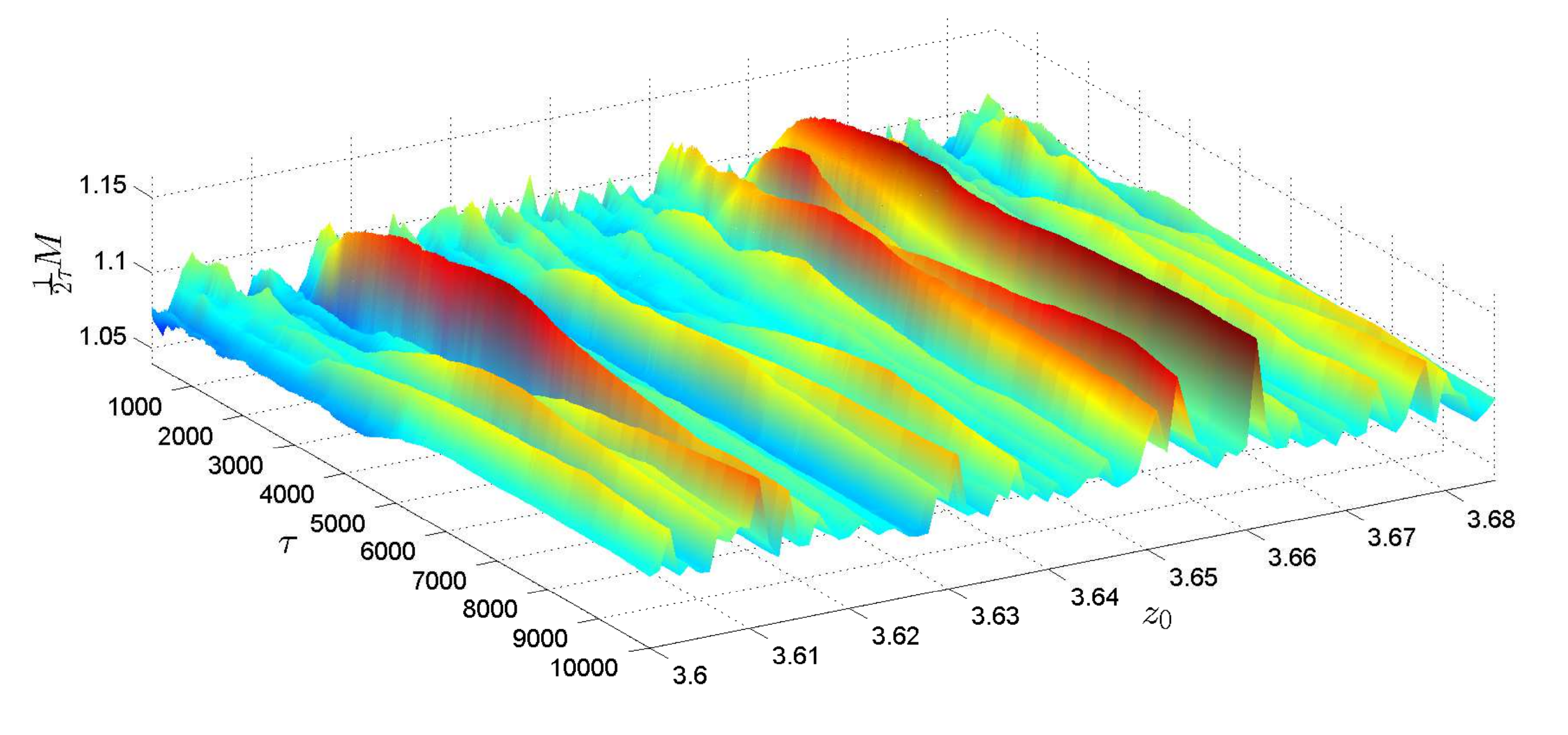}}
b){\includegraphics[scale = 0.42]{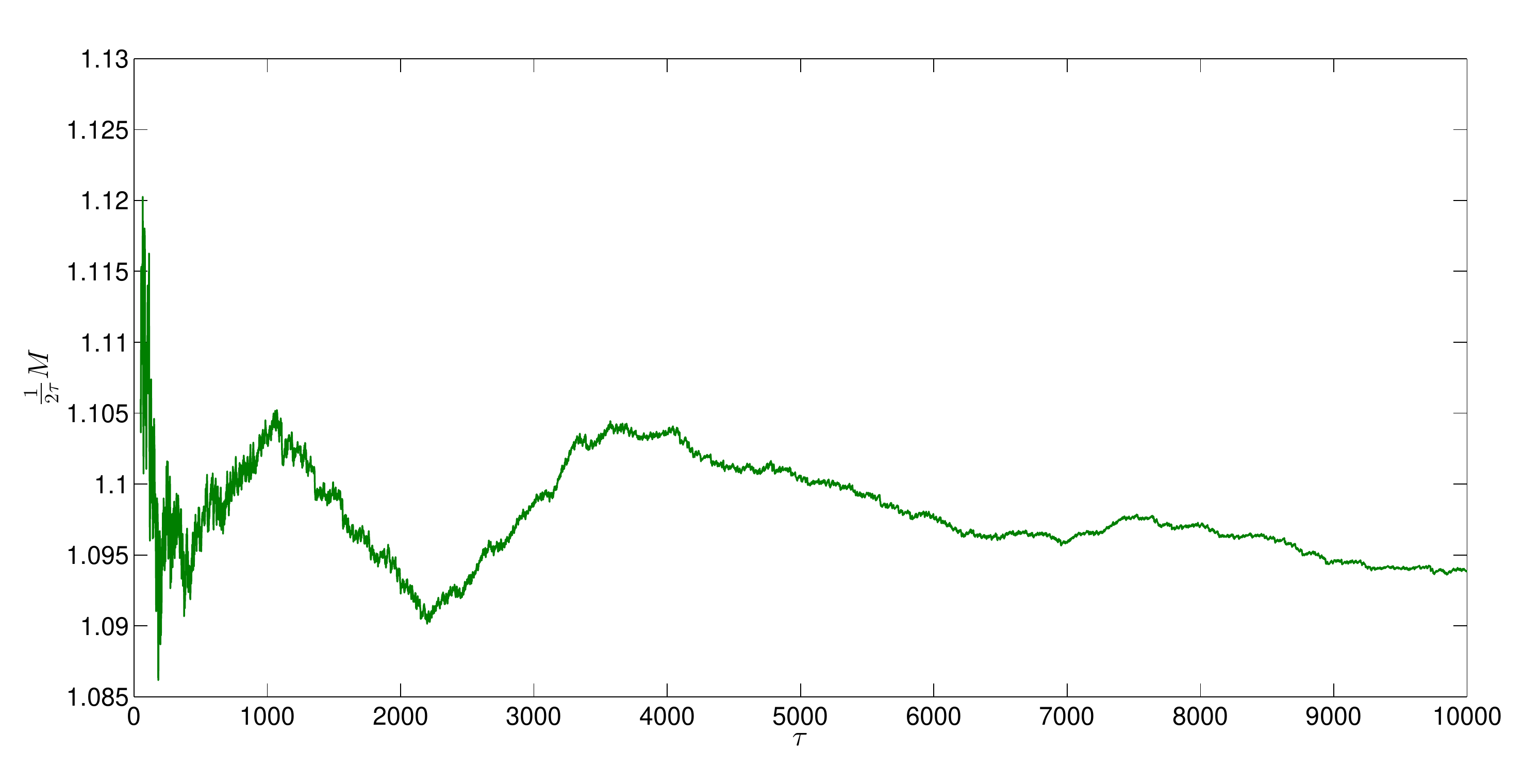}}
\caption{a) Time average evolution of $M$ in the range $\tau \in (0,10000]$ for a subset of the line  (\ref{ic_line}) in the chaotic region; b) Time average evolution of  green initial condition in the range $\tau \in (0,10000]$ .}
\label{t_avg_conv_chaos}
\end{figure}

\begin{figure}[htbp!]
\centering
\includegraphics[scale = 0.6]{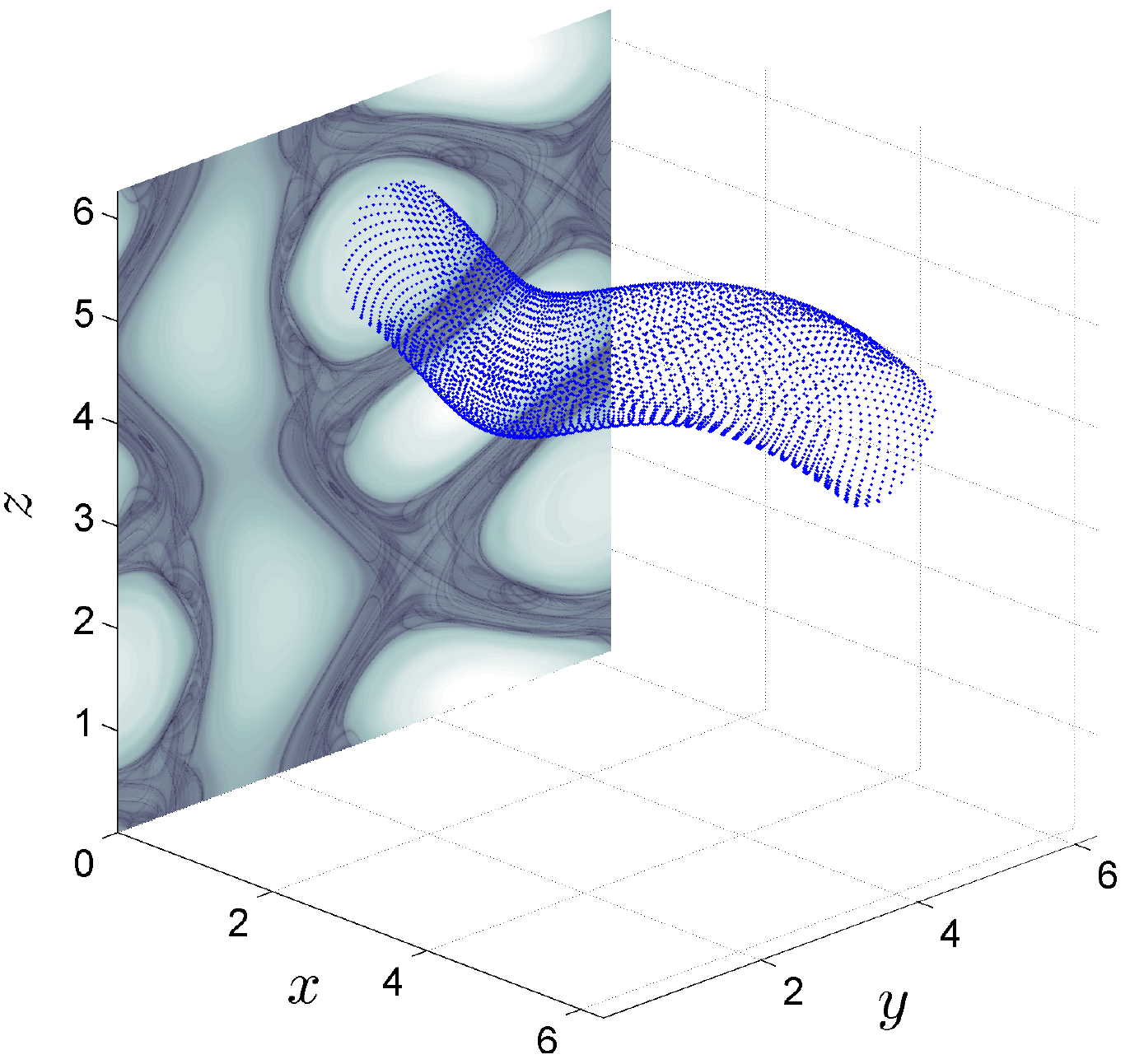}
\caption{Trajectory of an initial condition inside the elliptic region which displays the associated invariant set.}
\label{traj_inv_set}
\end{figure}

\section{Objectivity and phase space structure}
\label{sec:LD_obj}

The utility of LDs for revealing phase space structure has been questioned in the literature \citep{haller15} as a result of them not having the property of objectivity. Briefly, a scalar valued and time-dependent, function is said to be {\em objective} if it is invariant under Galilean coordinate transformations. In other words, the pointwise values of a function are the same at points that are transformed under a Galilean transformation, for each value of the time variable, see, e.g., \cite{Truesdell2004,haller2016}. Other accepted definitions for objectivity are given in the literature in terms of consistency between frames  \cite{jfm,pea}  but in this section we base our discussion on the objectivity definition given above as it is the one considered when debating LDs performance. Certainly in physics many scalar valued functions describing physical quantities, such as energy, or the magnitude of angular momentum, should be invariant under Galilean transformations. But this is {\em not} a property that is desirable for any tool designed to reveal phase space structure since phase space structure may  not  be invariant under Galilean coordinate transformations. We demonstrate this in the following example.

We consider the simplest possible dynamical system on the plane, the zero vector field:
\begin{equation}
\dot{{\bf x}}=0 \quad,\quad {\rm where} \quad {\bf x} \in {\rm \mathbb{R}^2} \label{sys1}
\end{equation}

This represents a system at rest. We apply a Galilean transformation to this vector field, i.e., a rotation 
$\underline{\bf x}=R(t)^{T}{\bf x}$, where $R(t)^T$ denotes  the transpose of the orthogonal matrix with angular speed $\omega=1$:
\begin{eqnarray}
R(t)= \left( \begin{array}{cc}
\cos  t & -\sin t  \\
\sin t & \cos t \end{array} \right)
\nonumber
\end{eqnarray}

\noindent In this rotating frame, \eqref{sys1} takes the form: 
\begin{equation}
\begin{cases}
\dot{\underline{x}} = \underline{y} \\
\dot{\underline{y}} = -\underline{x}
\end{cases}
\label{sys2}
\end{equation}

\noindent
The  phase portrait of  (\ref{sys1}) consists entirely of fixed points. The phase portrait of (\ref{sys2}), which represents a particular case of an harmonic oscillator,  consists of a one-parameter family of invariant circles (see \cite{arnold73} pp 44-45).  Clearly the phase space structure of these two dynamical systems is different, and as the LDs can be analytically computed for both dynamical systems, this fact is verified explicitly.

The arclength based LD, denoted by $M$,  measures the arclength of a trajectory through an initial condition in both forward and backward time, for a specified time.
For \eqref{sys1} this  is identically zero since  every point is a fixed point, and therefore the arclength of every trajectory is zero, regardless of the time for which it is computed. 
As for \eqref{sys2}, we recall results in Section \ref{sec:LD_elliptic} where it was shown that  $M=(2\tau)\rho$ for this example. The contours of $M$  are the same as the contours  of the Hamiltonian $H=\rho$, and therefore the contours of $M$  are in 1-1 correspondence with the trajectories of (\ref{sys2}). Therefore $M$ recovers the correct  phase space structure for both (\ref{sys1}) and (\ref{sys2}). Evidently, if $M$ were objective, i.e., the same for each of these vector fields, it would not recover the phase space structure for each vector field.}

It is instructive to consider what Lyapunov exponents, both finite and infinite time, would reveal for these examples. For \eqref{sys1} both the finite and infinite time Lyapunov exponents {\em for any} trajectory are zero. For \eqref{sys3} the infinite time Lyapunov exponents of every trajectory are zero and the value of the finite time Lyapunov exponents depend on the time interval over which they are computed. Therefore, we can make the following conclusions.

\begin{itemize}

\item While the infinite time Lyapunov exponents are objective, in the sense that they give the same values for \eqref{sys1} and \eqref{sys3}, they fail to reveal the phase space structure for \eqref{sys3}.

\item Finite time Lyapunov exponents are {\em not} objective. Their values, and hence  the phase space structure that they reveal, depend on the time interval over which they are computed.  A discussion of this can be found in \cite{brawigg, cnsns}.

\end{itemize}

We discuss next another example taken from \citep{haller2005,haller15,thesiswang} to show that LDs recover the correct phase space structure even when it is not evident in the instantaneous streamline phase portrait. We consider the following time-dependent dynamical system,
\begin{equation}
\dot{\mathbf x}= \left( \begin{array}{cc}
\sin 2 \omega t & \omega+\cos 2 \omega t  \\
-\omega+\cos 2\omega t & -\sin 2\omega t \end{array} \right){\mathbf x}
\label{sysob}
\end{equation}

\begin{figure}[htbp!]
\centering
a){\includegraphics[scale = 0.4]{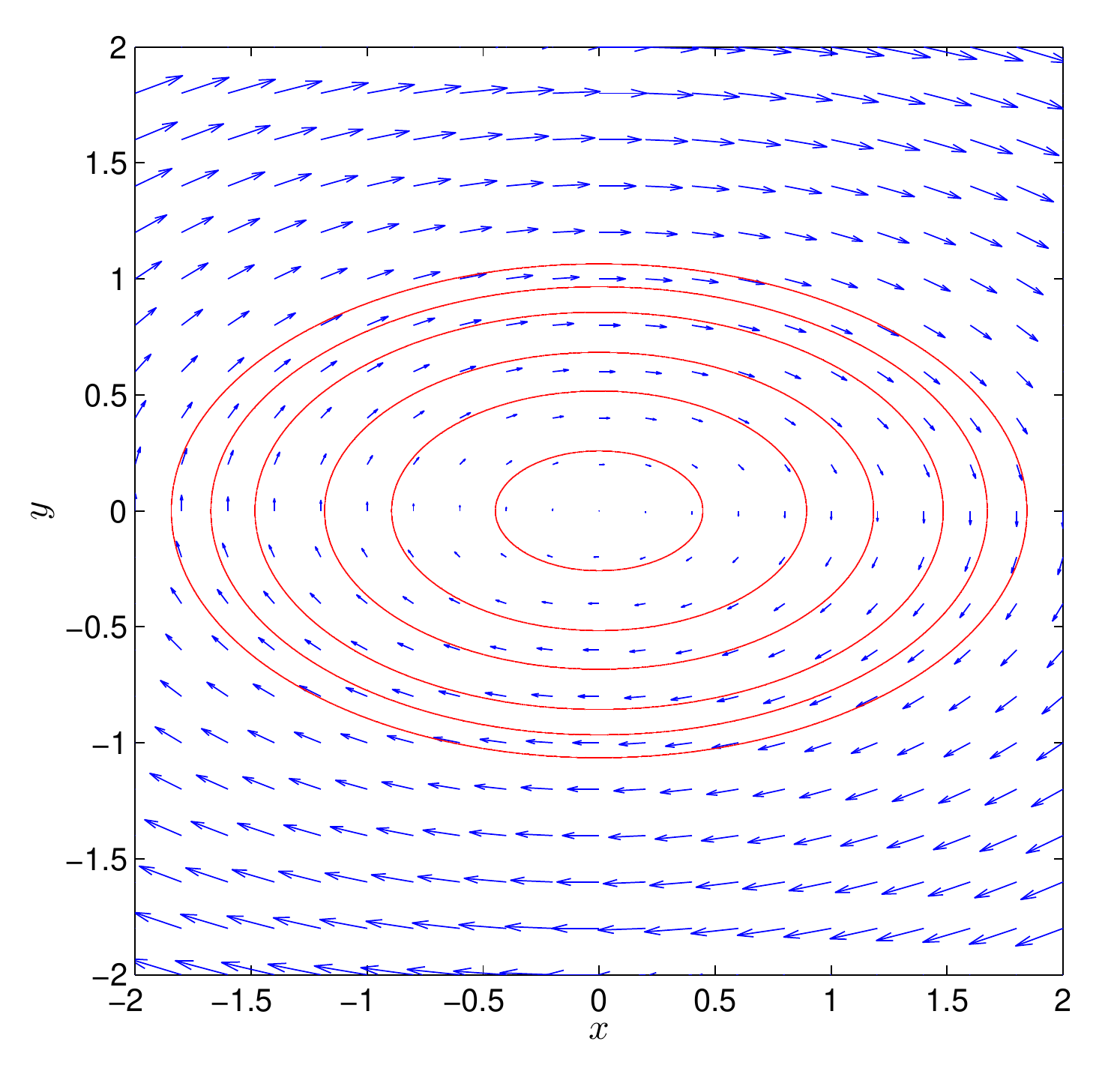}}
b){\includegraphics[scale = 0.4]{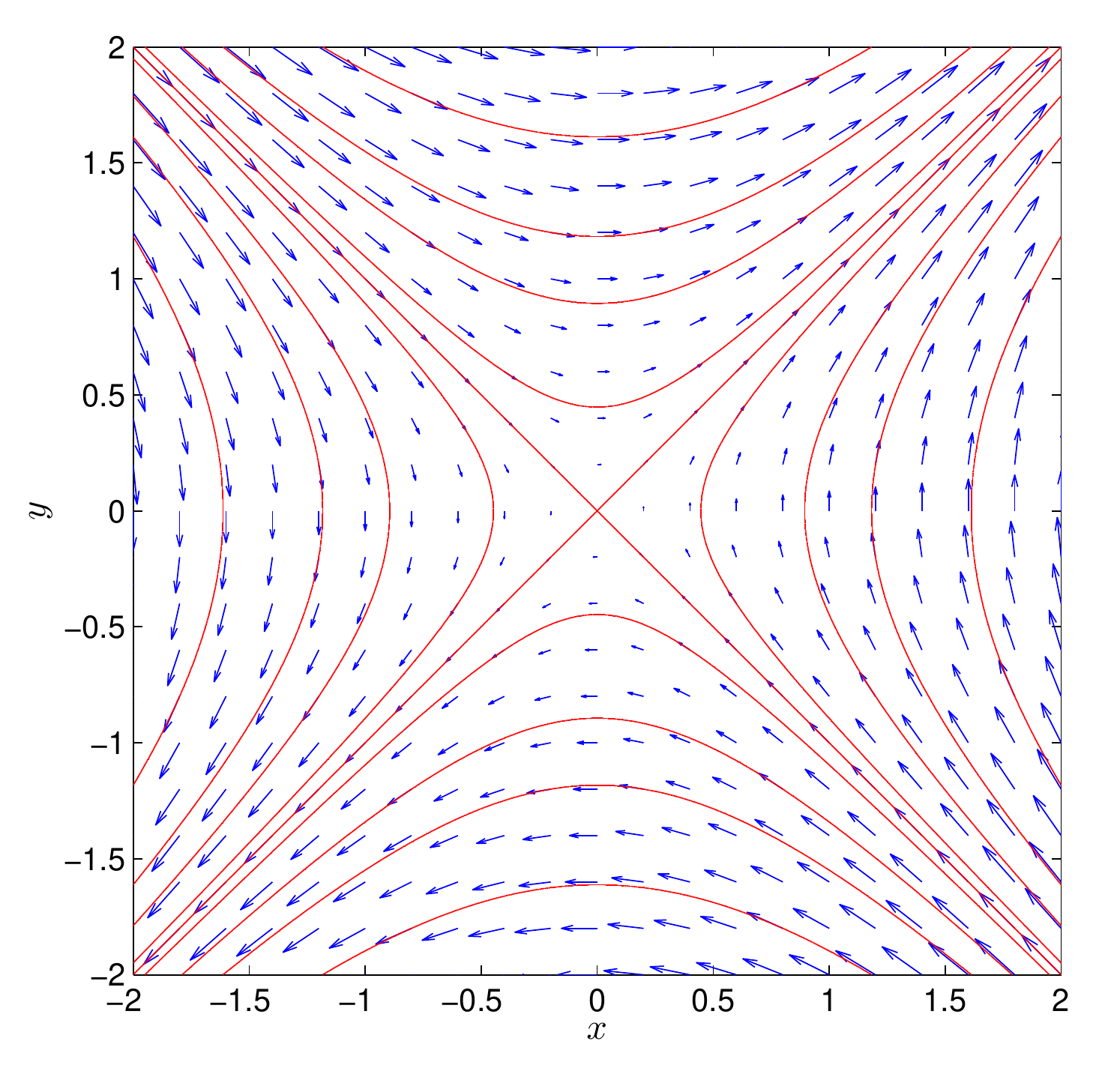}}
\caption{a) Streamlines and velocity field for system (\ref{sysob}) at $t=0$. b) Streamlines and velocity field for system (\ref{eq:lin_aut_saddle_point_ob}) at $t=0$.}
\label{obj}
\end{figure}

\noindent Figure \ref{obj}a) shows the instantaneous velocity fields and streamlines at $t=0$ for this example with $\omega=2$. They show a circulating pattern suggesting the presence of an eddy. However the following coordinate transformation,
\[\mathbf{x} = R(t)^T \underline{\mathbf{x}} \quad,\quad
R(t) = 
\begin{pmatrix}
\cos  \omega t & -\sin  \omega t  \\
\sin \omega  t & \cos  \omega t
\end{pmatrix}
\]

\noindent
where $R(t)^T$ is the transpose of $R(t)$, converts system (\ref{sysob}) into a stationary saddle,
\begin{equation}
\begin{cases}
   \dot{\underline{x}} = \underline{y}\\
   \dot{\underline{y}} = \underline{x} 
\end{cases}\;.
\label{eq:lin_aut_saddle_point_ob}
\end{equation}

\begin{figure}
\centering
{\includegraphics[scale = 0.35]{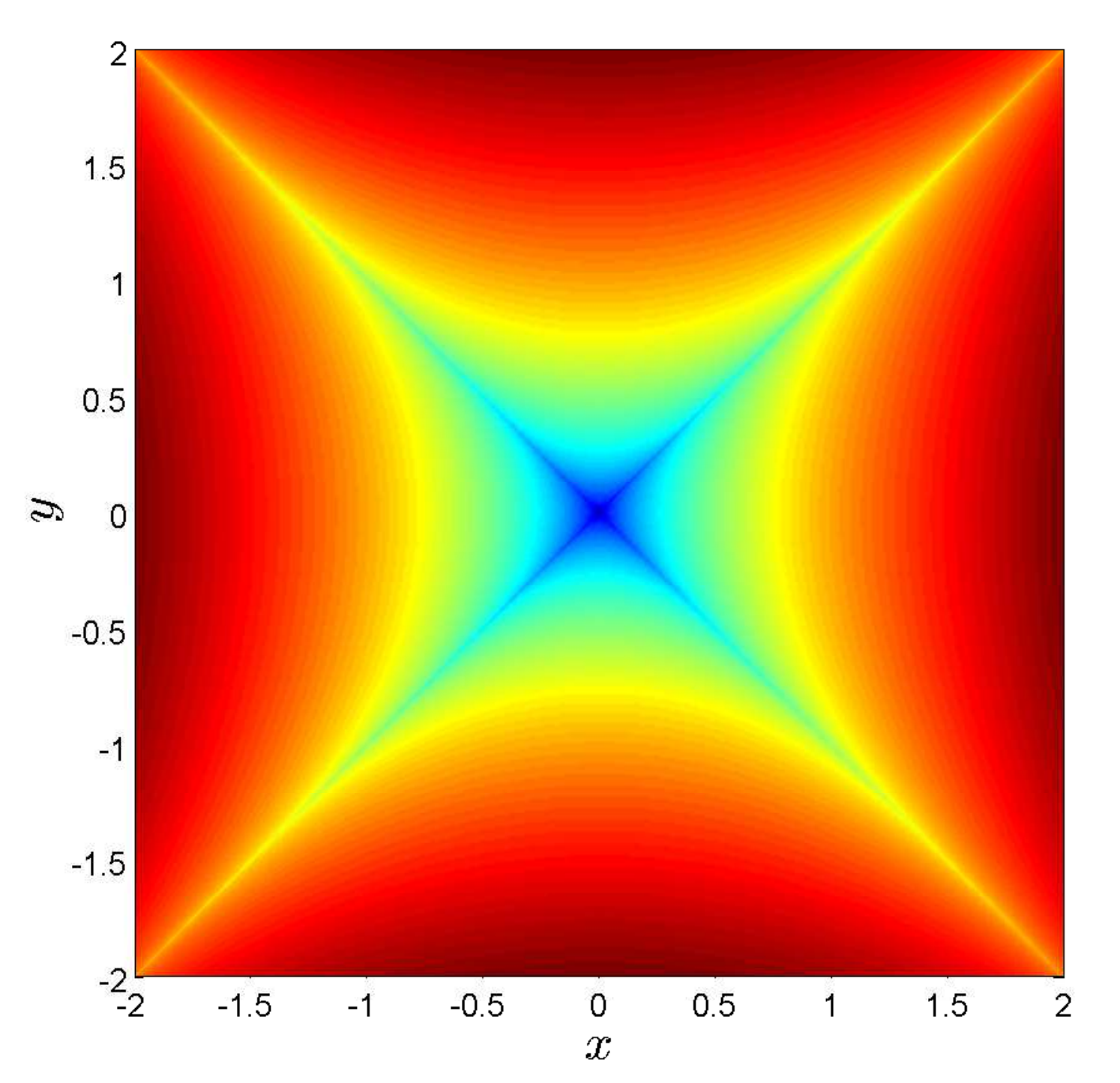}}
{\includegraphics[scale = 0.35]{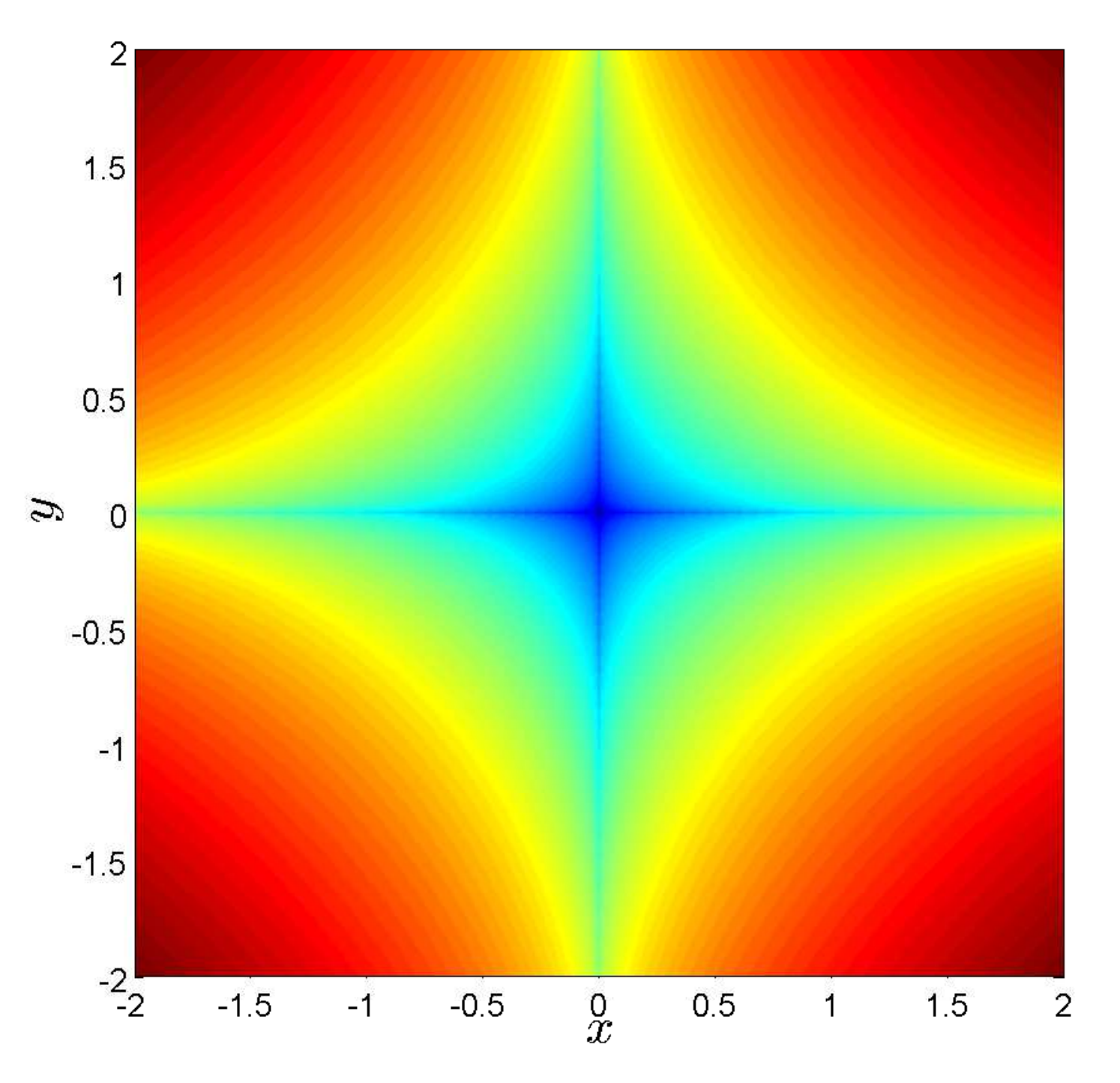}}
\caption{Contours of $M_{p=0.5}$ obtained for system (\ref{sysob}) with $\tau=10$ at succesive times a) $t=0$ and b) $t=\pi/8$.}
\label{obj2}
\end{figure}

Hence \eqref{sysob}, despite the structure revealed by the instantaneous streamline curves, is actually a rotating saddle point, i.e. a saddle point at the origin with rotating stable and unstable manifolds.  Precisely this structure is revealed by the LD $M_p$, as it is illustrated in   Fig.  \ref{obj2},   which shows contours of  $M_{p=0.5}$ at successive times $t=0, \pi/8$. The contours clearly reveal the rotating saddle point structure.

In this example $M_p$ is clearly not objective since the pointwise values of $M_p$ vary  from \eqref{sysob} to \eqref{eq:lin_aut_saddle_point_ob}.  Nevertheless, it is clear that {\em if} $M_p$ satisfied this criterion of objectivity, it would be the same for both systems and thus it would {\em not} distinguish between  the phase space structure for each of these very different  systems, providing an inconsistent description in both frames. 
Of course the values of $M_p$ at specific points of space certainly change with the reference frame, but the points at which $M_p$  is singular --which are the features containing the Lagrangian information-- are transformed with a smooth change of coordinates in the same manner in which the manifolds themselves are transformed. This was also illustrated for the rotating Duffing equation in \cite{jfm}.

Moreover, with respect to the question of the requirement of  objectivity in the context of techniques for revealing Lagrangian structure, it is instructive to note the following. 
Recently, in the context of fluid mechanics, a technique called Lagrangian-averaged vorticity deviation (LAVD) \citep{haller2016} has been developed. This technique,  by construction, has   the property of being invariant under Galilean transformations. Consequently, it does not distinguish between the phase portraits of \eqref{sys1} and \eqref{sys2}, as we now show.

The Lagrangian-averaged vorticity deviation is defined as follows,
\begin{equation}
{\rm LAVD}_{t_0}^{t_0+\tau}({\bf x}_0)=\int_{t_0}^{t_0+\tau} |w({\bf x}({\bf x_0},t_0,t),t)-\overline{w}(t)| dt
\end{equation}

\noindent
where the vorticity $w=\nabla \times {\bf v}$, is evaluated along the trajectory ${\bf x}({\bf x_0},t_0,t)$. In this expression $\overline{w}$
is the instantaneous spatial mean of the vorticity over $U(t)$:
$$
\overline{w}(t)=\frac{\int_{U(t)}   w({\bf x},t) dV}{{\rm vol}(U(t))}.
$$
where $U(t)$  is a domain invariant under the fluid flow and vol() denotes the volume. Let us evaluate LAVD for system  \eqref{sys2}. The vorticity is constant everywhere $w=-2$, in particular also along the trajectories. On the other hand let us consider the domain $U(t)=\{ (x,y) | \sqrt{x^2+y^2}\leq \rho_1\}$ which is invariant under the flow of  system  \eqref{sys2}. In this case it is easily found that:
$$
\overline{w}(t)=\frac{ -2 \int_{U(t)}   dV}{\pi  \rho_1^2}=-2.
$$
thus LAVD  is constantly zero on the whole domain for  \eqref{sys2}. This is also  clearly the case for \eqref{sys1} and thus  LAVD contour lines do not distinguish between these systems. 

\noindent

Finally, a reflection  on the property of objectivity --understood as a property of functions which preserve pointwise values under a Galilean transformation-- , and when it may be required, is useful.
From the point of view of distinguishing Lagrangian structures in velocity fields related by a Galilean coordinate transformation, our examples show that  objectivity is {\em not} a desirable property for a method to detect phase space structures in different frames.

\section{Conclusions}
\label{sec:conclusions}

This paper provides a theoretical framework for Lagrangian descriptors. In particular, the issues surrounding the ability of LDs to detect invariant stable and unstable manifolds of hyperbolic points are stated and clarified. This is accomplished by precisely defining the notion of "singular feature" and rigorously proving the presence  of  these features aligned with invariant manifolds in  four particular
cases: a hyperbolic saddle point for linear autonomous systems, a hyperbolic saddle point for nonlinear autonomous systems, a hyperbolic saddle point for linear nonautonomous systems and a hyperbolic saddle point for nonlinear nonautonomous systems. In order to achieve this goal we have
proposed a new way of constructing Lagrangian descriptors that keeps proofs simple.

We have also discussed well known  rigorous results  of the ergodic partition theory which are related to LDs. As a result we 
show the ability of LDs to highlight additional invariant sets, such as n-tori, by means of contour plots of converged averages. 

We have presented  the application of LDs to  the 3D ABC flow, in which it is shown how LDs locate simultaneously invariant manifolds, that are distinguishable as singular features, and invariant tori, visible from contour lines. The ability of LDs to highlight both manifolds and coherent eddy like or jet like structures had been noted in the literature \citep{amism11,wm14,ggmwm15}, although in this paper it is linked to previously known rigorous results.
We note however that these works dealt with aperiodic  flows and there is no generalization of the Birkhoff ergodic theorem to the case of aperiodically time-dependent vector fields, which is required for the ergodic partition theory.

Finally, we have provided a discussion of the topic of objectivity  in the context of Lagrangian descriptors. Specifically, we have analyzed their ability to provide the correct description of phase space structures
under Galilean transformation, as well as showing that the requirement of the objectivity property in general for tools in order to reveal phase space structures is not a desirable property. 

\section*{\bf Acknowledgements} The research of CL, FB-I, VJG-G and AMM is supported by the MINECO under grant 
MTM2014-56392-R. The research of SW is supported by  ONR Grant No.~N00014-01-1-0769.  We acknowledge support from 
MINECO: ICMAT Severo Ochoa project SEV-2011-0087.

\end{document}